\documentclass[12pt,a4paper,twoside]{article}
\usepackage{times}
\usepackage{amsmath,amssymb,amsbsy,amsgen,amsopn,amscd,amstext,amsxtra}
\usepackage{amsthm}
\usepackage{mathrsfs}
\usepackage{enumerate}
\usepackage{fancyhdr}
\newtheorem{Thm}{Theorem}[section]
\newtheorem{Prop}[Thm]{Propostion}
\newtheorem{Lem}[Thm]{Lemma}
\newtheorem{Def}[Thm]{Definition}
\newtheorem{Cor}[Thm]{Corollary}
\theoremstyle{remark}
\newtheorem*{Rem}{\textbf{Remark}}

\makeatletter
 
  \@addtoreset{equation}{section}
 \makeatother
\usepackage[dvipdfm]{color}
\usepackage[dvipdfm,bookmarks=true,linkbordercolor={0 1 1},
pdftitle=Article,
pdfsubject=Smoothing effects on real rank one symmetric spaces,
pdfauthor=Koichi Kaizuka,
pdfkeywords=Harmonic Analysis On Symmetric Spaces]{hyperref}
\usepackage[abbrev]{amsrefs}
\makeatletter
\def\address#1#2{\begingroup
\noindent\parbox[t]{7.8cm}{%
\small{\scshape\ignorespaces#1}\par\vskip1ex
\noindent\small{\itshape E-mail address}%
\/:#2\par\vskip4ex}\hfill%
\endgroup}%
\makeatother
\begin{document}
\title{\textbf{Smoothing effects of dispersive equations\\
               on real rank one symmetric spaces}}
\author{Koichi Kaizuka}
\date{}
\maketitle
%%%%%%%%%%%%%%%%%%%%%%%%%%%%%%%%%%
%         foot note              %
%%%%%%%%%%%%%%%%%%%%%%%%%%%%%%%%%%
\footnote[0]{%2000 MSC numbers
}
\footnote[0]{%key words and phrase
}
\footnote[0]{%Thanks
$^{*}$The author is supported by the JSPS Research
Fellowships for Young Scientists and the JSPS Grant-in-Aid 
for Scientific Research No.20-6882
}
%%%%%%%%%%%%%%%%%%%%%%%%%%%%%%%%%%%%%%%%%%%%%%%%%%%%%%%%%%%%%%%%%%%%%%
%  Abstract                                                          %
%%%%%%%%%%%%%%%%%%%%%%%%%%%%%%%%%%%%%%%%%%%%%%%%%%%%%%%%%%%%%%%%%%%%%%
\vspace*{-20pt}
\begin{abstract}
In this article we prove time-global smoothing effects of 
dispersive pseudodifferential equations with 
constant coefficient radially symmetric symbols on 
real rank one symmetric spaces of noncompact type. 
We also discuss gain of
regularities according to decay rates of 
initial values for the Schr\"{o}dinger evolution equation. 
We introduce some isometric operators and 
reduce the arguments to the well-known Euclidean case.
In our proof, Helgason's Fourier transform and the Radon transform
as an elliptic Fourier integral operator play crucial roles.
\end{abstract}
%%%%%%%%%%%%%%%%%%%%%%%%%%%%%%%%%%%%%%%%%%%%%%%%%%%%%%%%%%%%%%%%%%%%%%
%  Section 1 : Introduction                                          %
%%%%%%%%%%%%%%%%%%%%%%%%%%%%%%%%%%%%%%%%%%%%%%%%%%%%%%%%%%%%%%%%%%%%%%
\section{Introduction}
Let $X$ be a real rank one symmetric space of noncompact type.
In this paper, we consider the initial value problem
for dispersive equations of the form
%%% IVP on X %%%%
\begin{alignat}{2}
D_{t}u-a(D_{x})u
&=f(t,x)
&\quad\text{in}\quad
&\mathbb{R}\times X, 
 \label{eq:IVPX1}\\
u(0,x)
&=\varphi (x) 
&\quad\text{in}\quad 
&X, 
 \label{eq:IVPX2}
\end{alignat}
where $u(t,x)$ is a complex-valued unknown function of
$(t,x)\in\mathbb{R}\times X$, $f(t,x)$ and $\varphi(x)$ are given
functions, $i=\sqrt{-1}$, $\partial_{t}=\partial / \partial t$,
$D_{t}=-i\partial_{t}$, and $a(D_{x})$ is a pseudodifferential operator
with a real-valued symbol on $X$
(the precise definition shall be given later).
In case that $a(D_{x})=-\Delta_{X}$,
where $\Delta_{X}$ is the Laplace-Beltrami operator on $X$,
the equation (\ref{eq:IVPX1}) becomes the
Schr\"{o}dinger evolution equation.

The purpose of this paper is to study the smoothness of the solution
to (\ref{eq:IVPX1})-(\ref{eq:IVPX2}) with additional assumptions on
the symbol of the pseudodifferential operator $a(D_{x})$.

Roughly speaking, by introducing some isometric operators, 
we can translate all the time-global smoothing estimates 
on the one-dimensional Euclidean space into those 
on the real rank one symmetric space $X$.

We see the Radon transform $\mathcal{R}$ as an elliptic 
Fourier integral  operators, and recover the regularity 
of a function $u$ on $X$ from $\mathcal{R}u$ and 
the canonical relation of $\mathcal{R}^{*}$.
Then we can also prove 
the gain of regularities for the solutions 
of the Scr\"{o}dinger evolution equations 
on the real rank one symmetric space $X$ from those 
on the one-dimensional Euclidean space.

First we will review some known results 
for dispersive equations on Euclidean spaces, 
mainly time-global spatially-local smoothing effects
for real-principal-type pseudodifferential equations.
We consider the initial value problem of the form
\begin{alignat}{2}
D_{t}u-a(D_{x})u
&= F(t,x)
&  \quad\text{in}\quad
&  \mathbb{R}\times \mathbb{R}^{n}, 
   \label{eq:IVPE1}\\ 
u(0,x)
&= \phi(x)
&  \quad\text{in}\quad
&  \mathbb{R}^{n}, \label{eq:IVPE2}
\end{alignat}
where $u(t,x)$ is a complex-valued unknown function 
of $(t,x)\in \mathbb{R}\times \mathbb{R}^{n}$,
$F(t,x)$ and $\phi (x)$ are given functions, 
$a(\xi)$ is a real-valued function at most
polynomial growth at infinity belongs to
$C^{1}(\mathbb{R}^{n})$
with $\nabla_{\xi} a(\xi)\neq 0$ for any
$\xi\neq 0$, and the operator $a(D_{x})$
is defined by
\begin{equation*}
a(D_{x})v(x)
=(2\pi)^{-n}
 \int_{\mathbb{R}^{2n}}
  e^{i(x-y)\xi}
  a(\xi)v(y)
 dyd\xi 
\end{equation*}
for an appropriate function $v$
on $\mathbb{R}^{n}$.

Since $a(\xi)$ is real-valued, 
the initial value problem 
(\ref{eq:IVPE1})-(\ref{eq:IVPE2}) is 
$L^{2}$-well-posed, that is, for any 
$\phi\in L^{2}(\mathbb{R}^{n})$ 
and for any 
$F\in L_{\mathrm{loc}}^{1}(\mathbb{R};L^{2}(\mathbb{R}^{n}))$, 
(\ref{eq:IVPE1})-(\ref{eq:IVPE2})
possesses a unique solution 
$u\in C(\mathbb{R};L^{2}(\mathbb{R}^{n}))$.

In the case $a(\xi)=|\xi|^{2}$,
i.e.\ $a(D_{x})=-\Delta_{\mathbb{R}^{n}}$,
the corresponding equation is the Schr\"{o}dinger evolution equation
and there exist many related works for its dispersive properties.

One of the origin of the studies for dispersive properties 
is the study for well-posedness of the initial value problem
of the KdV type equations by Kato ~\cite{Kato} around 1983.

In ~\cite{KatoYajima} Kato and Yajima obtained
the following estimate in case $n\geq 3$:  
\begin{equation*} 
 \left\|
      \langle x\rangle^{-1}\langle D_{x}\rangle^{1/2}
      e^{-it\Delta_{\mathbb{R}^{n}}}\phi
 \right\|_{L^{2}(\mathbb{R}\times\mathbb{R}^{n})}
\leq
 C \|\phi\|_{L^{2}(\mathbb{R}^{n})}.
\end{equation*}
Here we put 
$\langle x\rangle  
 =(1+|x|^{2})^{1/2}$ 
and 
$\langle D_{x}\rangle
=\mathcal{F}_{\mathbb{R}^{n}}^{-1} 
 \langle\xi\rangle
 \mathcal{F}_{\mathbb{R}^{n}}$,
$|D_{x}|
=\mathcal{F}_{\mathbb{R}^{n}}^{-1} 
 |\xi|
 \mathcal{F}_{\mathbb{R}^{n}}$.

In ~\cite{BenKla} Ben-Artzi and Klainerman 
also showed the following type estimate
in case $n\geq 3$:
\begin{equation*}
 \left\|
     \langle x\rangle^{-\delta}
     |D_{x}|^{1/2}
     e^{-it\Delta_{\mathbb{R}^{n}}}\phi
 \right\|_{L^{2}(\mathbb{R}\times\mathbb{R}^{n})}
\leq
 C_{\delta} 
 \|\phi\|_{L^{2}(\mathbb{R}^{n})},
\end{equation*}
where $\delta>1/2$ and $C_{\delta}>0$ 
is a positive constant depending on $\delta$. 

In ~\cite{Sugimoto} Sugimoto obtained another type of 
time-global smoothing estimates for generalized Schr\"{o}dinger 
operator not only for the homogeneous solutions 
but also for the inhomogeneous solutions.

Later in ~\cite{Hoshiro} Hoshiro obtained time-local smoothing
effects for general real-principal polynomial symbols, that is,
$
  a(\xi)
=\sum_{|\alpha|\leq m}
  a_{\alpha}\xi^{\alpha}
$
$
(m\geq 2,a_{\alpha}\in\mathbb{R})
$
with its principal symbol
$
 a_{m}(\xi)
=\sum_{|\alpha|=m}
  a_{\alpha}\xi^{\alpha}
$
satisfies the dispersive condition $\nabla_{\xi} a_{m}\neq 0$
for any $\xi\neq 0$. 
Especially, in ~\cite{Hoshiro} Hoshiro
proved that the dispersive condition
$\nabla_{\xi} a_{m}(\xi)\neq 0\,(\xi\neq 0)$ 
is a necessary condition for the time-local 
spatially-local smoothing effects to hold.

Here we remark that the dispersive condition 
corresponding to a non-trapping condition 
for a Hamilton orbit generated 
by the principal symbol, i.e.\
$$
\nabla_{\xi}a_{m}(\xi)
\neq 0 \quad (\xi\neq 0)
\,\Leftrightarrow\,
|x+t\nabla_{\xi}a_{m}(\xi)|
\to
\infty
\,\,
(|t|\to\infty)
\quad\text{for\,any}\,
(x,\xi)
\in
T^{*}\mathbb{R}^{n}\setminus 0 .
$$

Recently, in ~\cite{Chihara} Chihara obtained 
time-global smoothing estimates for
real-principal-type positively homogeneous symbols of
degree $m>1$, that is, 
real-valued functions
$a(\xi) 
 \in C^{\infty}(\mathbb{R}^{n}\setminus \{0\})
 \cap C^{1}(\mathbb{R}^{n})$
such that
$\nabla a(\xi)\neq 0$,
$
 a(\xi)
=|\xi|^{m} 
 a\left( \xi / |\xi| \right)
$
$(\xi \neq 0)$.
Set 
$p(\xi)=|\xi|^{(m-1)/2}$
if $n\geq 2$,
$p(\xi)=a'(\xi)|\xi|^{-(m-1)/2}$
if $n=1$.
Then the estimate is as follows:
$$
 \left\|
     \langle x\rangle^{-\delta}
     p(D_{x})u
 \right\|_{L^{2}(\mathbb{R}\times\mathbb{R}^{n})}
\leq
C_{\delta}
\left(
     \|\phi\|_{L^{2}(\mathbb{R}^{n})}
     +\left\|
          \langle x\rangle^{\delta}
          |D_{x}|^{(m-1)/2}
          F
      \right\|_{L^{2}(\mathbb{R}\times\mathbb{R}^{n})}
\right),
$$
where $\delta >1/2$ and $C_{\delta}$ is a positive
constant depending on $\delta$.

%%%%% Notes for smoothing effects %%%%%%%%%%%%%%%%%%%%%%%%%%%%%%%
Above results are all for equations with constant coefficients.
The methods of the proofs are mainly based on some
Fourier restriction theorem in the homogeneous case and
the limiting absorption principle in the inhomogeneous case.
Also, for equations with variable coefficients, there are
so many results related to smoothing effects.
For example, in ~\cite{Doi1} Doi proved that the non-trapping condition is
necessary for the spatially-local smoothing effect of 
Schr\"{o}dinger evolution groups on complete Riemannian manifolds.
(Also see ~\cite{CKS}, ~\cite{Doi2}, ~\cite{Doi3}, ~\cite{Doi0}, ~\cite{Chihara}.) 

Now we focus our interest on time-global smoothing effects for the 
Schr\"{o}dinger evolution equation. 
Then above time-global estimates are essentially 
flat-Euclidean case.
For other noncompact complete Riemannian manifolds 
time-global smoothing effects had not been obtained
except for a special case.
In ~\cite{RodTao} Rodnianski and Tao obtained
a time-global smoothing estimate for the 
Schr\"{o}dinger evolution equation on $\mathbb{R}^{3}$
with compact metric perturbations.
From the geometrical point of view for the Schr\"{o}dinger evolution 
equation, if we assume some ``nice'' geometric structures 
for the Riemannian manifold and the non-trapping condition 
for the geodesic flow,
time-global smoothing effects may be expected. 
Moreover, to investigate the relation 
between geometrical conditions and the time-global smoothing 
effects is one of the interesting and deep problem.

%%%%% Notes for Harmonic analysis %%%%%%%%%%%%%%%%%%%%%%

In this paper Helgason's Fourier transform plays crucial roles
as in the Euclidean case.
See section \ref{PRE} for the notation and 
the rigourous definition for harmonic analysis on symmetric spaces.

In ~\cite{Helgason1}, ~\cite{Helgason2} and ~\cite{Helgason3} 
Helgason introduced the Fourier transform on symmetric spaces, 
and also proved the Plancherel formula, the inversion formula, 
and the Paley-Wiener theorem for his Fourier transform. 
After his pioneering works, harmonic analysis on symmetric spaces have been
studied actively and applied to various fields of
mathematics by many people.

%%%%%%%%%  Paragraph of  Main results  %%%%%%%%%%%%%%%%%%%%%%%%%%%%%
On the symmetric space $X$ of noncompact type, 
Helgason's Fourier transform is defined by  
$$
 \mathcal{F}u(\lambda ,b)
=\int_{X}
  e^{(-i\lambda +\rho)(A(x,b))}
  u(x)
 dx
\quad 
\text{for}
\,\,
u\in\mathscr{D}(X)
$$
for $(\lambda,b)\in\mathfrak{a}^{*}\times B$.
And the Fourier transform is invertible by
the following formula: 
$$
u(x)
=\iint_{\mathfrak{a}_{+}^{*}\times B}
  e^{(i\lambda +\rho)(A(x,b))}
  \mathcal{F}u(\lambda,b)
  |\boldsymbol{c}(\lambda)|^{-2}
 d\lambda db
\quad 
\text{for\,any}
\,\,x\in X.
$$
In the rest of this section, we will state our main results.

%%%%%%%% Main results %%%%%%%%%%%%%%%%%%%%%%%%%%%%%%%%%%%%%%
For any real-valued function
$a(\lambda)
 \in
 C(\overline{\mathfrak{a}^{*}_{+}})
 \cap
 C^{1}(\mathfrak{a}^{*}_{+})$
at most polynomial order at infinity,
the pseudodifferential operator $a(D_{x})$
is defined by
$$
 a(D_{x})v(x)
=\iint_{\mathfrak{a}^{*}_{+}\times B}
   e^{(i\lambda+\rho)(A(x,b))}
   a(\lambda)
   \mathcal{F}v(\lambda ,b)
   |\boldsymbol{c}(\lambda)|^{-2}
 d\lambda db.
$$
for an appropriate function $v$ on $X$.

Since $a(\lambda)$ is real-valued, 
the initial value problem 
(\ref{eq:IVPX1})-(\ref{eq:IVPX2}) is 
$L^{2}$-well-posed, that is, for any 
$\varphi\in L^{2}(X,dx)$ 
(abbreviated $L^{2}(X)$ in the sequel)
and for any 
$f\in L_{\mathrm{loc}}^{1}(\mathbb{R};L^{2}(X))$, 
(\ref{eq:IVPX1})-(\ref{eq:IVPX2}) possesses a 
unique solution 
$u\in C(\mathbb{R};L^{2}(X))$.
Moreover, the unique solution $u$
is explicitly given by 
\begin{align*}
 u(t,x)
&=e^{ita(D_{x})}\varphi (x)+iGf(t,x),  \\
 e^{ita(D_{x})}\varphi (x)
&=\iint_{\mathfrak{a}_{+}^{*}\times B}
   e^{(i\lambda +\rho)(A(x,b))}
   e^{ita(\lambda)}
   \mathcal{F}\varphi (\lambda ,b)
   |\boldsymbol{c}(\lambda)|^{-2}
  d\lambda db,  \\
 Gf(t,x)
&=\int_{0}^{t}
   e^{i(t-\tau)a(D_{x})}
   f(\tau ,x)
  d\tau . 
\end{align*}

Here we introduce a ``weighted'' $L^{2}$-space
on $X$ as follows.

%%% Definition of weighted L2-space %%%%%%%%%%%%%%%%
\begin{Def}
 \label{Def:WL2space}
Let $\delta\in\mathbb{R}$, and for all 
$u,v\in\mathscr{D}(X)$, define an inner metric
\begin{equation*}
(u,v)_{L^{2,\delta}(X)}
=w^{-1}
 \iint_{\mathfrak{a}^{*}\times B}
   \langle D_{\lambda}\rangle^{\delta}
   \big( 
      \mathcal{F}u(\lambda ,b)
      \boldsymbol{c}^{-1}(\lambda) 
   \big) 
  \overline
  { \langle D_{\lambda}\rangle^{\delta}
    \big( 
       \mathcal{F}v(\lambda ,b)
       \boldsymbol{c}^{-1}(\lambda) 
    \big) }
 d\lambda db, 
\end{equation*}
and ``weighted'' $L^{2}$-norm
$
 \|u\|_{L^{2,\delta}(X)}
=(u,u)^{1/2}_{L^{2,\delta}(X)}$.
Let $L^{2,\delta}(X)$ be the completion 
of the pre-Hilbert space
$\big( 
   \mathscr{D}(X),
   (\cdot ,\cdot )_{L^{2,\delta}(X)}   
 \big)
$
(well-definedness will be discussed in section \ref{WL2}).
\end{Def}

In what follows we assume that $X$ is real rank one.
From Tits' classification ~\cite{Tits} of
isotoropic homogeneous manifolds,
the real rank one symmetric spaces of noncompact type
are exactly the  following four types of
noncompact complete Riemannian manifolds,
$$
H^{n}(\mathbb{R}),
H^{n}(\mathbb{C}),
H^{n}(\mathbb{H}),
\text{or}\, 
H^{2}(\mathbb{O}),
$$ 
those are, the real, complex, quaternion hyperbolic space or 
the octave hyperbolic plane respectively, 
where $\mathbb{H}$ is the quaternion
and $\mathbb{O}$ is the octave or called the Cayley algebra. 

Here we state our first main results.
%%%%%%%%% Theorem of Smoothing Effects %%%%%%%%%%%%%%%%%%%%%%%%
\begin{Thm}
\label{Thm:SmEf}
$\mathrm{(i)}$ Let $a(\lambda)\in C^{1}(\mathfrak{a}^{*})$ 
             and $p(\lambda)\in C^{0}(\mathfrak{a}^{*})$ be real-valued
             functions at most polynomial order at infinity.
             Suppose that the pseudodifferential operator
             $a(D_{H})$ on $\mathfrak{a}$ causes 
             a time-global estimate for homogeneous solutions, 
             that is, 
             there exist positive constants 
             $\delta >1/2$ and $C_{\delta}$ such that 
             $$
             \left\| 
                 \langle H\rangle^{-\delta}
                 p(D_{H})e^{ita(D_{H})}
                 \phi
             \right\|_{L^{2}(\mathbb{R}\times\mathfrak{a},dtdH)}
             \leq C_{\delta}
                  \|\phi\|_{L^{2}(\mathfrak{a},dH)}.
             $$
             Then we have
             $$
             \left\| 
                 p(D_{x})
                 e^{ita(D_{x})}
                 \varphi
             \right\|_{L^{2}(\mathbb{R};L^{2,-\delta}(X))}
             \leq w^{1/2}
                  C_{\delta}
                  \|\varphi\|_{L^{2}(X)}.
             $$
$\mathrm{(ii)}$ Let $a(\lambda)\in C^{1}(\mathfrak{a}^{*})$
             and $q(\lambda)\in C^{0}(\mathfrak{a}^{*})$ be real-valued
             functions at most polynomial order at infinity. 
             Suppose that the pseudodifferential operator $a(D_{H})$
             on $\mathfrak{a}$ causes a time-global estimate
             for inhomogeneous solutions, that is,
             there exist positive constants $\delta >1/2$,
             $C_{\delta}$ and a cut off function
             $\chi\in C^{\infty}(\mathfrak{a}^{*})$ with
             $\chi(\lambda )=0\,(|\lambda |\leq 1)$,
             $\chi(\lambda )=1\,(|\lambda |\geq 2)$ such that 
             \begin{align*}
             \left\| 
                 \langle H\rangle^{-\delta}
                 \chi (D_{H})q(D_{H})
                 \int_{0}^{t}
                 e^{i(t-\tau )a(D_{H})}
                 F(\tau)
                 d\tau
             \right\|_{L^{2}(\mathbb{R}\times\mathfrak{a},dtdH)}
             \leq C_{\delta}
                  \|\langle H\rangle^{\delta}F
                  \|_{L^{2}(\mathbb{R}\times\mathfrak{a},dtdH)},
              \end{align*}
             Then there exists a positive constant $C_{\delta , \chi}$
             such that
             \begin{align*}
              \left\| 
                  \chi(D_{x})q(D_{x})
                  \int_{0}^{t}
                  e^{i(t-\tau)a(D_{x})}
                  f(\tau)
                  d\tau
              \right\|_{L^{2}(\mathbb{R};L^{2,-\delta}(X))}
             \leq C_{\delta ,\chi}
                  \| f\|_{L^{2}(\mathbb{R};L^{2,\delta}(X))}.
             \end{align*}
 
 \end{Thm}
%%%%% Corolally for Smoothing effects %%%%%%%%%%%%%%%%%%%%%%%%%%%%%
\begin{Cor}
 \label{Cor:SmEf}
Let $a(D_{x})$ be a polynomial of the Laplace-Beltrami operator
$\Delta_{X}$ of real coefficients,
then for any $\delta >1/2$ and
$\chi\in C^{\infty}(\mathfrak{a}^{*})$
as in Theorem \ref{Thm:SmEf}
there exist positive constants
$C_{\delta},C_{\delta, \chi}$
such that
\begin{align*}
 \left\| 
     |a'(D_{x})|^{1/2}
     e^{ita(D_{x})}
     \varphi
 \right\|_{L^{2}(\mathbb{R};L^{2,-\delta}(X))}
&\leq C_{\delta}
       \| \varphi\|_{L^{2}(X)}, \\
 \left\| 
     \chi(D_{x})a'(D_{x})
     \int_{0}^{t}
       e^{i(t-\tau)a(D_{x})}
       f(\tau)
     d\tau
 \right\|_{L^{2}(\mathbb{R};L^{2,-\delta}(X))}
 &\leq C_{\delta ,\chi}
       \| f\|_{L^{2}(\mathbb{R};L^{2,\delta}(X))}.
 \end{align*}
 \end{Cor}
%%%%%%%%%% Remark for Smoothing effects %%%%%%%%%%%%%%%%%
\begin{Rem}
  Needless to say, the estimates above are meaningless without
  some regularity estimates for the ``weighted'' $L^{2}$-space
  $L^{2,\delta}(X)$ and we will study local regularities and
  basic properties of $L^{2,\delta}(X)$ in Section \ref{WL2}.
  \end{Rem}

Here we prepare standard notation for micorolocal analysis. 
%%%%%%%%% Definition of wave front set %%%%%%%%%%%%%%%%%%%
\begin{Def}
Let $M$ be an $n$-dimentional smooth manifold and
$s$ be any real number.
A distribution $v\in\mathscr{D}'(M)$ is in the local Sobolev space
$H^{s}_{\mathrm{loc}}(M)$
if and only if
$\rho^{*}(\chi v)\in H^{s}(\mathbb{R}^{n})$
for any local chart $(\rho, U)$ and
cut-off function $\chi\in C_{0}^{\infty}(U)$.
Also the distribution $v\in\mathscr{D}'(M)$
is in $H^{s}$ microlocally near
$(x_{0},\xi_{0})\in T^{*}M\setminus 0$
if and only if
there exists a local chart $(\rho ,U)$,
cut-off funtion $\chi\in C_{0}^{\infty}(U)$
with $\chi(x_{0})\neq 0$ and
a conic neighborhood
$V\subset\mathbb{R}^{n}\setminus \{0\}$ of
$\eta_{0}=(\rho^{-1})^{*}\xi_{0}$
such that
$\langle \eta \rangle^{s}
 \mathcal{F}_{\mathbb{R}^{n}}
 (\rho^{*}(\chi v))(\eta)
 \in L^{2}(V,d\eta)$.
Let $\operatorname{WF}^{s}(v)$
denote the elements of $T^{*}M\setminus 0$
such that $v$ is not in $H^{s}$ microlocally.
$\operatorname{WF}^{s}(v)$ is called the Sobolev
wave front set of $v\in\mathscr{D}'(M)$ of
order $s\in\mathbb{R}$.
\end{Def}

Also, we have the gain of regularity for
the Schr\"{o}dinger evolution equation as following.
%%%%%%%%% Theorem for Gain of regurality %%%%%%%%%%%%%%%%%
\begin{Thm}
 \label{Thm:GRe}
Suppose $\varphi\in L^{2}(X)$ satisfies
the following condition for some positive integer $k$ and
open set $\Theta\subset B$:
\begin{equation}
 \label{eq:AsGR}
\left\| 
    \langle D_{\lambda}\rangle^{k}
    \left\{ 
        \mathcal{F}\varphi(\lambda,b)
        \boldsymbol{c}(\lambda)^{-1}
    \right\}
\right\|_{L^{2}(\mathfrak{a}^{*}\times \Theta ,d\lambda db)}
< \infty .
\end{equation}
Then for all $x\in X$, and $b\in \Theta_{x}$ we have
\begin{alignat*}{2}
 (x;\pm\omega_{b}(x))
&\notin 
 \operatorname{WF}^{k+1/2}
 (e^{-it\Delta}\varphi)
&\quad
 \text{for a.e.}
 \quad
&t\neq 0,\\
 (x;\pm\omega_{b}(x))
&\notin 
 \operatorname{WF}^{k}
 (e^{-it\Delta}\varphi)
&\quad\text{for any}\quad
&t\neq 0.
\end{alignat*}
In particular, if $\phi\in L^{2,k}(X)$, then we have 
\begin{alignat*}{2}
 e^{-it\Delta_{X}}\varphi
&\in 
 H_{\mathrm{loc}}^{k+1/2}(X)
&\quad\text{for a.e.}\quad
&t\neq 0,\\
e^{-it\Delta_{X}}\varphi
&\in 
 H_{\mathrm{loc}}^{k}(X)
&\quad\text{for any}\quad
&t\neq 0.
\end{alignat*}
Here
$
 \omega_{b}(x)
=d_{x}\big( A(x,b)\big)
 \in\Lambda^{1}(X)
$, 
$
 \Theta_{g\cdot o}
=\Theta 
  \cap 
  \big( 
     g
     \circ 
     (-\mathrm{id}_{B})
     \circ 
     g^{-1}\Theta
  \big)
$,
$g\in G$.
\end{Thm}
%%%%%%% Theorem Gain Continuity %%%%%%%%%%%%%%%%%%%%%%%%%%%%%%%%%%%
\begin{Thm}
\label{Thm:GainConti}
For any $k\in\mathbb{N}$ 
and $\delta > 1/2$, 
we have the following continuous maps: 
\begin{alignat*}{2}
 L^{2,k}(X)\ni \varphi 
&\mapsto
 t^{k}
 \langle D_{x}\rangle^{k+1/2}
 e^{-it\Delta_{X}}\varphi
&\hspace{4pt}
 \in
&\hspace{4pt}
 L^{2}_{\text{loc}}
 (\mathbb{R};
 L^{2,-k-\delta}(X)),\\
 L^{2,k}(X)\ni \varphi 
&\mapsto
 t^{k}
 \langle D_{x}\rangle^{k}
 e^{-it\Delta_{X}}\varphi
&\hspace{4pt} 
 \in
&\hspace{4pt}
 C(\mathbb{R};
 L^{2,-k}(X)).
\end{alignat*}
\end{Thm}

\begin{Rem}
\begin{itemize}
\item[$(\mathrm{i})$] 
The results in 
Theorem \ref{Thm:GRe} 
and \ref{Thm:GainConti}
may be contained in Doi's results
in ~\cite{Doi2},
but our proof takes different approach 
and is rather simple.
So we will treat those as theorems.
\item[$(\mathrm{ii})$]
Here we remark that the condition (\ref{eq:AsGR}) 
corresponds to a decaying along some family of geodesics.
In fact, we can rewrite the norm by using 
a pseudodifferential operator $\Lambda$ 
and the Radon transform $\mathcal{R}$ as following
($\Lambda$ and $\mathcal{R}$ are given in the next section):
$$ 
\left\| 
    \langle D_{\lambda}\rangle^{k}
    \left\{ 
        \mathcal{F}\varphi(\lambda,b)
        \boldsymbol{c}(\lambda)^{-1}
    \right\}
\right\|_{L^{2}(\mathfrak{a}^{*}\times\Theta ,d\lambda db)}
=\left\| 
     \langle H\rangle^{k}
     \Lambda\mathcal{R}\varphi
 \right\|_{L^{2}(\mathfrak{a}\times\Theta ,d\xi(H,b))}.
$$
\item[($\mathrm{iii}$)]
We can restate the results in Theorem 
\ref{Thm:GRe} by using geodesics:
For a maximal geodesic $\gamma :\mathbb{R}\to X$, 
if $b_{+}=\gamma (+\infty)$,
$b_{-}=\gamma (-\infty) \in \Theta$, 
then for any $s\in\mathbb{R}$ we have
\begin{alignat*}{2}
 (\gamma (s);
  \omega_{b_{\pm}}(\gamma(s)))
&\notin 
 \operatorname{WF}^{k+1/2}
 (e^{-it\Delta}\varphi)
&\quad
 \text{for a.e.}
 \quad
&t\neq 0,\\
 (\gamma (s);
  \omega_{b_{\pm}}(\gamma(s)))
&\notin 
 \operatorname{WF}^{k}
 (e^{-it\Delta}\varphi)
&\quad\text{for any}\quad
&t\neq 0.
\end{alignat*}
\end{itemize}
\end{Rem}

%%%%%% contents of this article %%%%%%%%%%%%%%%
Finally, this paper is organized as follows.
In Section \ref{PRE} 
we review Helgason's harmonic analysis on 
symmetric spaces.
In Section \ref{WL2} 
we prove local regularites for elements of 
weighted $L^{2}$-space.
In Section \ref{TGSE} 
we prove Theorem \ref{Thm:SmEf}.
In Section \ref{GOR} 
we prove Theorems \ref{Thm:GRe} and \ref{Thm:GainConti}.
%%%%%%%%%  
\vspace{1mm}
%%%%%%%%%%%%%%%%%%%%%%%%%%%%%%%%%%
%%%%%%%  Acknowledgments   %%%%%%%
%%%%%%%%%%%%%%%%%%%%%%%%%%%%%%%%%%
\section*{Acknowledgments} 
The author would like to express his sincere
gratitude to Hiroyuki Chihara
for his great encouragement, valuable and precise advices.
He would like to thank the members of the seminar, 
Ryuichiro Mizuhara and Eiji Onodera
for their kindly advices and helpful discussions.
%%%%%%%%%%%%%%%%%%%%%%%%%%%%%%%%%%%%%%%%%%%%%%%%%%%%%%%%%%%%%%%%%%%%%%
%  Section 2 : Preliminaries                                         %
%%%%%%%%%%%%%%%%%%%%%%%%%%%%%%%%%%%%%%%%%%%%%%%%%%%%%%%%%%%%%%%%%%%%%%
\section{Preliminaries\label{PRE}}
In this section we introduce Helgason's harmonic analysis 
and prepare some lemmas needed later. 
One can consult with Helgason's books 
~\cite{HeGGA}, ~\cite{HeGASS} and ~\cite{HeDLS}
on the basic facts 
on harmonic analysis on symmetric space.
Also, we will establish some tools on symmetric spaces, 
related to the proof of the theorems.

In this section, 
$X$ is an arbitrary rank symmetric space
unless we impose some conditions in addition.

\begin{figure}[h]
\begin{center}
\unitlength 0.1in
\begin{picture}( 30.9600, 25.0600)( 45.5700,-56.5300)
% DOT 1 0 3 0
% 4 6400 4400 6400 4400 6400 4400 6400 4400
% 
\special{pn 13}%
\special{sh 1}%
\special{ar 6400 4400 10 10 0  6.28318530717959E+0000}%
\special{sh 1}%
\special{ar 6400 4400 10 10 0  6.28318530717959E+0000}%
\special{sh 1}%
\special{ar 6400 4400 10 10 0  6.28318530717959E+0000}%
\special{sh 1}%
\special{ar 6400 4400 10 10 0  6.28318530717959E+0000}%
% CIRCLE 2 0 3 0
% 4 6400 4400 7286 3514 7286 3514 7286 3514
% 
\special{pn 8}%
\special{ar 6400 4400 1254 1254  0.0000000 6.2831853}%
% CIRCLE 2 0 3 0
% 4 6990 3810 7286 3514 7286 3514 7286 3514
% 
\special{pn 8}%
\special{ar 6990 3810 420 420  0.0000000 6.2831853}%
% LINE 2 2 3 0
% 2 5514 5286 7286 3514
% 
\special{pn 8}%
\special{pa 5514 5286}%
\special{pa 7286 3514}%
\special{dt 0.045}%
% DOT 1 0 3 0
% 3 6694 3514 7286 3514 7286 3514
% 
\special{pn 13}%
\special{sh 1}%
\special{ar 6694 3514 10 10 0  6.28318530717959E+0000}%
\special{sh 1}%
\special{ar 7286 3514 10 10 0  6.28318530717959E+0000}%
\special{sh 1}%
\special{ar 7286 3514 10 10 0  6.28318530717959E+0000}%
% STR 2 0 3 0
% 3 6289 4326 6289 4400 5 0
% $o$
\put(62.8900,-44.0000){\makebox(0,0){$o$}}%
% STR 2 0 3 0
% 3 5367 3293 5367 3367 5 0
% $B=\mathbb{S}^{1}$
\put(53.6700,-33.6700){\makebox(0,0){$B=\mathbb{S}^{1}$}}%
% STR 2 0 3 0
% 3 5800 3974 5800 4050 5 0
% $X=H^{2}(\mathbb{R})$
\put(58.0000,-40.5000){\makebox(0,0){$X=H^{2}(\mathbb{R})$}}%
% STR 2 0 3 0
% 3 6593 3340 6593 3412 5 0
% $x$
\put(65.9300,-34.1200){\makebox(0,0){$x$}}%
% STR 2 0 3 0
% 3 7541 3291 7541 3365 5 0
% $b=kM$
\put(75.4100,-33.6500){\makebox(0,0){$b=kM$}}%
% STR 2 0 3 0
% 3 6260 3686 6260 3760 5 0
% $\xi (x,b)$
\put(62.6000,-37.6000){\makebox(0,0){$\xi (x,b)$}}%
% STR 2 0 3 0
% 3 7023 4411 7023 4485 5 0
% $k\exp (A(x,b))\!\cdot\! o$
\put(70.2300,-44.8500){\makebox(0,0){$k\exp (A(x,b))\!\cdot\! o$}}%
% DOT 1 0 3 0
% 2 6694 4106 6694 4106
% 
\special{pn 13}%
\special{sh 1}%
\special{ar 6694 4106 10 10 0  6.28318530717959E+0000}%
\special{sh 1}%
\special{ar 6694 4106 10 10 0  6.28318530717959E+0000}%
% VECTOR 2 0 3 0
% 2 6822 4400 6704 4138
% 
\special{pn 8}%
\special{pa 6822 4400}%
\special{pa 6704 4138}%
\special{fp}%
\special{sh 1}%
\special{pa 6704 4138}%
\special{pa 6714 4208}%
\special{pa 6726 4188}%
\special{pa 6750 4192}%
\special{pa 6704 4138}%
\special{fp}%
\end{picture}%
 \caption{The horocycle on the Poincar\'{e} disc}
 \label{Fig:Horocycle}
\end{center}
\end{figure}
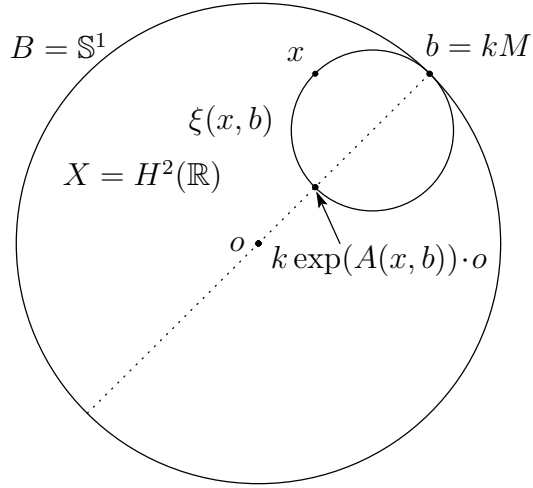

 Let $X$ be a Riemannian symmetric space of noncompact type, 
that is, $X=G/K$ where $G$ is a noncompact, connected, 
semisimple Lie group with finite center and
$K$ is a maximal compact subgroup of $G$. 
Let $\mathfrak{g}$ and $\mathfrak{k}$ be the Lie algebra of 
$G$ and $K$ respectively. 
Let $B$ be the Killing form on $\mathfrak{g}$, that is, 
$B(X,Y)=\operatorname{Trace}
         (\operatorname{ad} (X)
          \operatorname{ad} (Y))$
for $X,Y\in \mathfrak{g}$,
where $\operatorname{ad}(X)Y=[X,Y]$, 
and $\theta$ be the Cartan involution associated 
with the Cartan decomposition 
$\mathfrak{g}=\mathfrak{k}+\mathfrak{p}$. 
Let $\langle X,Y\rangle=-B(X,\theta{Y})$,
then $\langle\cdot, \cdot\rangle$ defines
an inner metric on $\mathfrak{g}$.
Let $\mathfrak{a}\subset\mathfrak{p}$ be a maximal abelian 
subspace and $\mathfrak{a}^{*}$ its dual. 
Then $l=\dim\mathfrak{a}$ is
called the real rank of $G$,  or $X$. 
For $\alpha\in\mathfrak{a}^{*}$, put 
$$
\mathfrak{g}_{\alpha}
=\{X\in\mathfrak{g};[H,X]=\alpha (H)X,H\in\mathfrak{a}\}
$$
and $m_{\alpha}=\dim\mathfrak{g}_{\alpha}$.
If $\alpha\neq 0$ and $\mathfrak{g}_{\alpha}\neq\{0\}$, 
then $\alpha$ is called a restricted root and denotes all 
restricted roots by $\Sigma$. 
Let $\mathfrak{g}_{\mathbb{C}}$ and
$\mathfrak{a}^{*}_{\mathbb{C}}$ denote the complexification 
of $\mathfrak{g}$ and $\mathfrak{a}^{*}$, respectively. 
If $\lambda , \mu\in\mathfrak{a}^{*}_{\mathbb{C}}$, 
let $H_{\lambda}\in\mathfrak{a}_{\mathbb{C}}$ be determined
by $\lambda (H)=\langle H_{\lambda},H \rangle$
for $H\in\alpha$ and put
$\langle \lambda,\mu\rangle
 =\langle H_{\lambda},H_{\mu}\rangle$. 
Since $B$ is positive definite on $\mathfrak{p}$,
we put
$|\lambda|=\langle\lambda,\lambda \rangle^{1/2}$
for 
$\lambda\in\mathfrak{a}^{*}$
and
$|X|=\langle X,X\rangle^{1/2}$
for 
$X\in\mathfrak{p}$. 
Then the natural identification
$(D\pi_{G})_{e}|_{\mathfrak{p}}:
 \mathfrak{p}\rightarrow T_{o}X$
(denote $o=eK$)
induces the left $G$-invariant Riemannian metric on $X$.
Let $dx$ be the Riemannian measure, 
and $\Delta_{X}$ be the Laplace-Beltrami operator on $X=G/K$.
Let $\mathfrak{a}'$ be the open subset of $\mathfrak{a}$  
where all the restricted roots are $\neq 0$. 
Fix a Weyl chamber $\mathfrak{a}^{+}$ in 
$\mathfrak{a}'$ that is a connected component of 
$\mathfrak{a}'$, and call $\mathfrak{a}\in\Sigma$  
is positive if it is positive on $\mathfrak{a}^{+}$. 
Let $\Sigma^{+}$ denotes the set of all positive roots. 
Put 
$\rho
=(1/2)\sum_{\alpha\in\Sigma^{+}}
      m_{\alpha}\alpha$, 
$\mathfrak{n}
=\sum_{\alpha\in\Sigma^{+}}
 \mathfrak{g}_{\alpha}$
and let $N$ denote the
corresponding analytic subgroup of $G$. 
Then we obtain an Iwasawa decomposition $G=KAN$ of $G$
and 
$\mathfrak{g}
=\mathfrak{k}
+\mathfrak{a}
+\mathfrak{n}
$
of $\mathfrak{g}$.
Each $g\in G$ can be uniquely written
$g=\kappa (g)\exp (H(g))n(g)$,
by $\kappa (g)\in K$, $H(g)\in\mathfrak{a}$,
and $n(g)\in N$.
Let $M$ denote the centralizer of $A$ in $K$,
$M'$ the normalizer of $A$ in $K$,
and the $W$ factor group $M'/M$, called the Weyl group. 
The group $W$ acts as a group of linear transformations 
of $\mathfrak{a}^{*}$ by 
$(s\lambda)(H)=\lambda(s^{-1}\cdot H)$ 
for $H\in\mathfrak{a}$, $\lambda\in\mathfrak{a}^{*}$  
and $s\in W$, where $g\cdot X=\operatorname{Ad} (g)X$ for 
$g\in G$, $X\in\mathfrak{g}$. 
Let $w$ denote the order of $W$. 
We fix an orthonormal basis on $\mathfrak{a}$ 
with respect to the Killing form and its dual basis
on $\mathfrak{a}$, then we can regard 
$\mathfrak{a}$ and $\mathfrak{a}^{*}$ as the 
Euclidean spaces of dimension $l$ respectively. 
The Killing form induces the Euclidean measures on 
$\mathfrak{a}$ and $\mathfrak{a}^{*}$, multiplying 
these by the factor $(2\pi)^{-l/2}$ 
we obtain the invariant measures $dH$ and $d\lambda$ 
on $\mathfrak{a}$ and $\mathfrak{a}^{*}$ respectively.
Put $B=K/M=G/MAN$. 
Let $\mathfrak{m}$ be the Lie algebra of $M$ and 
$\mathfrak{l}$ be the orthogonal complement of 
$\mathfrak{m}$ in $\mathfrak{k}$ with respect to 
the Killing form. 
Since $\langle\cdot,\cdot\rangle$ is strictly positive on 
$\mathfrak{l}$, the natural identification
$(D\pi_{K})_{e}|_{\mathfrak{l}}:
 \mathfrak{l}\rightarrow T_{b_{0}}B$
(denote $b_{0}=eM$) 
induces the left $K$-invariant 
Riemannian metric on $B$. 
Let $db$ be the left $K$-invariant measure on $B$ 
normalized so that total mass equals one. 
Put $\Xi=G/MN$, called the horocycle space of $X$,
and let $d\xi$ be the invariant measure on $\Xi$.
We can naturally identify $\Xi$ with
$\mathfrak{a}\times K/M$ by the diffeomorphism
$\mathfrak{a}\times K/M\ni (H,kM)\mapsto
 k\cdot\exp(H)MN\in\Xi$,
then we can write $d\xi =e^{2\rho (H)}dHdb$.
For $x\in X, b\in B$ let $\xi(x,b)$ be the horocycle
passing through the point $x=g\cdot o\in X$ with normal 
$b=kM\in B$, and let 
$A(x,b)=-H(g^{-1}k)\in\mathfrak{a}$ be the composite distance
from the origin to $\xi (x,b)$.
(See Figure \ref{Fig:Horocycle}.)

%%%%%%%%%%%%%%% Figure of horocycle %%%%%%%%%%%%%%%%%%%%%%
%%%%%%%%%%%%%%%%%%%%%%%%%%%%%%%%%%%%%%%%%%%%%%%%%%%%%%%%

%%%%% The Helgason Fourier transform %%%%%%%%%%%%%%%%%%%%%%%
Now, we have the following fundamental identity
$$
 \Delta_{X}
 (e^{(-i\lambda+\rho)(A(x,b))})
=-(|\lambda|^{2}+|\rho|^{2})
  e^{(-i\lambda+\rho)(A(x,b))},
$$ 
then the Helgason Fourier transform of
$f\in\mathscr{D}(X)$ is given by
$$
\mathcal{F}f(\lambda,b)
=\int_{X}
  e^{(-i\lambda+\rho)(A(x,b))}f(x)
 dx,
 \quad
 (\lambda,b)\in\mathfrak{a }_{\mathbb{C}}^{*}\times B.
$$

Helgason proved the inversion formula and the Plancherel theorem
as follows:
%%%%%%%%%%%%%% Theorem for the Fourier transform %%%%%%%%%%%%%%%%%%%%%%%%%%%%%%%%%%%%%
\begin{Thm}
 \label{Thm:Helgason}
We have the followings:
\begin{itemize}
\item \textbf{Inversion formula} ; 
 For $f\in\mathscr{D}(X)$, the Fourier transform
 is inverted by
 $$
 f(x)=w^{-1}
      \iint_{\mathfrak{a}^{*}\times B}
       e^{(i\lambda+\rho)(A(x,b))}
       \mathcal{F}f(\lambda ,b)
       |\boldsymbol{c}(\lambda)|^{-2}
      d\lambda db
 \quad \text{for\, all}\; x\in X , 
 $$
 where $\boldsymbol{c}$ is Harish-Chandra's
 $\boldsymbol{c}$-function.
\item \textbf{The Plancherel theorem} ;
 The Fourier transform is extended to the unitary
 isomorphism
 $$
 \mathcal{F}:
 L^{2}(X,dx)
 \rightarrow
 L_{W}^{2}(\mathfrak{a}^{*}\times B,
           w^{-1}|\boldsymbol{c}(\lambda)|^{-2}
           d\lambda db ) ,
 $$
 where
 $ L_{W}^{2}(\mathfrak{a}^{*}\times B,
            |\boldsymbol{c}(\lambda)|^{-2}
            d\lambda db )$
 consisting of all elements of
 $L^{2}(\mathfrak{a}^{*}\times B,
        |\boldsymbol{c}(\lambda)|^{-2}
        d\lambda db )$
 satisfying the following condition
 \begin{equation}\label{eq:Winv}
  \int_{B}
       e^{(is\lambda +\rho)(A(x,b))}
       \psi(s\lambda ,b)
      db
 =\int_{B}
       e^{(i\lambda +\rho)(A(x,b))}
       \psi(\lambda ,b)
      db
 \end{equation}
 for all $s\in W$, and a.e. $x\in X$, $b\in B$.
 \end{itemize}
\end{Thm}
%%%%%%%%% Remark for the Fourier transform %%%%%%%%%%%%%%
\begin{Rem}
 By the $W$-invariantness of the Fourier images, we can
 rewrite the inversion formula as follows:
 $$
 f(x)=\iint_{\mathfrak{a}^{*}_{+}\times B}
        e^{(i\lambda+\rho)(A(x,b))}
        \mathcal{F}f(\lambda ,b)
        |\boldsymbol{c}(\lambda)|^{-2}
      d\lambda db.
 $$
 Also we have the following unitary isomorphism:
 $$
 \mathcal{F}:
 L^{2}(X,dx)
 \rightarrow
 L^{2}(\mathfrak{a}^{*}_{+}\times B,
           |\boldsymbol{c}(\lambda)|^{-2}
           d\lambda db ) .
 $$
\end{Rem}
Next, we introduce the Schwartz space $\mathscr{S}(X)$ defined
by Eguchi and Okamoto in ~\cite{EgOk}.
They characterized the image of it by the Fourier transform.
Here we identify the functions on $X=G/K$ with
the right $K$-invariant functions on $G$.
For a function $f$ on $X$, and for $X,X'\in\mathfrak{g}$
we put 
\begin{align*}
f(X ;\, g)&=\left. \frac{d}{dt}
                f(\exp (tX)g) \right|_{t=0}, \\
f(g\,;X')        &=\left. \frac{d}{dt}
                f(g\exp (tX)) \right|_{t=0}.
\end{align*}
And let $U(\mathfrak{g}_{\mathbb{C}})$ be the universal
enveloping algebra of $\mathfrak{g}_{\mathbb{C}}$,
then for $X$ and $X'\in U(\mathfrak{g}_{\mathbb{C}})$,
$f(X; g)$ and $f(g;X')$ are naturally defined respectively 
as the homomorphic extension of the above definition.

For $\lambda\in \mathfrak{a}_{\mathbb{C}}^{*}$,
let $\varphi_{\lambda}(x)\in \mathscr{E}(X)$ be
the elementary spherical function, that is,
$$
\varphi_{\lambda}(x)
=\int_{B}e^{(i\lambda +\rho)(A(x,b))}db,
$$
then $\varphi_{\lambda}(x)$ is a $K$-invariant eigenfunction on X.
For $\varphi_{0}(x)$, following basic properties are known.
(For the detail, see ~\cite{GaVa}, ~\cite{HaCh}.)  
%%%%% Proposition of Elementary SphericalFunction %%%%%%%%%%%%%%
\begin{Prop}
 \label{Prop:ProESF}
We have 
\begin{itemize}
\item[$(\mathrm{i})$]
  There exist a positive constant $C_{1}$,
 and a positive integer $N_{1}\in\mathbb{N}$
 such that
 \begin{equation*}
 e^{-\rho (\log a)}
 \leq \varphi_{0}(a)
 \leq C_{1}e^{-\rho (\log a)}(1+\sigma (a))^{N_{1}}
 \label{eq:EsESF1}
 \end{equation*}
 for any $a\in A^{+}=\exp \mathfrak{a}^{+}$.
\item[$(\mathrm{ii})$]
 There exists a positive integer $N_{2}\in\mathbb{N}$
 such that
 \begin{equation*}
 \int_{G}\varphi_{0}(g)^{2}(1+\sigma (g))^{-N_{2}}dg
 < \infty .
 \end{equation*}
\end{itemize}
\end{Prop}
%%%%%%%%%%

Now we define the Schwartz space on the symmetric space $X$.
%%%%%%%%%% Definition of Schwartz space on X %%%%% %%%%%%%%%
\begin{Def}
 \label{Def:SSS}
Let $\mathscr{S} (X)$ denote the set of $C^{\infty}$
functions $f$ on $X=G/K$ satisfying 
$$
\tau_{N,X,X'}(f)
=\sup_{g\in G}
  \left\{
         |f(X;g;X')|\varphi_{0}(g)^{-1}
         (1+\sigma (g))^{N}
  \right\}
< \infty ,
$$
for any $N\in\mathbb{Z}_{\geq 0}$, and
$X,X'\in U(\mathfrak{g}_{\mathbb{C}})$.
Then $\mathscr{S}(X)$ becomes a Fr\'{e}chet space with
the seminorms $\tau_{N,X,X'}$
$(N\in\mathbb{Z}_{\geq 0}, X,X'\in U(\mathfrak{g}_{\mathbb{C}}))$.
We can easily see that
$$
\mathscr{D}(X)\hookrightarrow
\mathscr{S}(X)\hookrightarrow
L^{2}(X),
$$
where all inclusions are continuous and dense embeddings.
\end{Def}
%%%%%%%%%

We also define the Schwartz space on the ``phase space''
$\mathfrak{a}^{*}\times B$.

%%%%%%%%% Definition of Schwartz space on phase space  %%%%%%%%%%%%
\begin{Def}
 \label{Def:SSF}
Let $\mathscr{I}(\mathfrak{a}\times B)$
denote the set of $C^{\infty}$ functions $\psi$
on $\mathfrak{a}\times B$ which satisfying 
$$
\nu_{N,\alpha,m}
=\sup_{\mathfrak{a}^{*}\times B}
  \left\{
    \left|
        (\partial_{\lambda}^{\alpha}
         \Delta_{B}^{m}
         \psi)
        (\lambda ,b)
    \right|
    (1+|\lambda |)^{N}
  \right\}
< \infty ,
$$
for any
$\alpha\in\mathbb{Z}_{\geq 0}^{l}$, $m,N\in\mathbb{Z}_{\geq 0}$.
Let $\mathscr{I}_{W}(\mathfrak{a}^{*}\times B)$
denote all the $W$-invariant elements $\psi$ in
$\mathscr{I}(\mathfrak{a}^{*}\times B)$, i.e.\
$\psi \in \mathscr{I}(\mathfrak{a}^{*}\times B)$
and satisfying 
$$
 \int_{B}e^{(is\lambda+\rho )(A(x,b))}
        \psi(s\lambda ,b)
 db
=\int_{B}e^{(i\lambda+\rho )(A(x,b))}
        \psi(\lambda ,b)
 db
$$
for all $s\in W$, $\lambda\in \mathfrak{a}^{*}$, and $x\in X$.
Then both $\mathscr{I}(\mathfrak{a}^{*}\times B)$ and
$\mathscr{I}_{W}(\mathfrak{a}^{*}\times B)$ are Fr\'{e}chet
spaces with the seminorms $\nu_{N,\alpha,m}$
$(\alpha\in\mathbb{Z}_{\geq 0}^{l}$, $m,N\in\mathbb{Z}_{\geq 0})$.
\end{Def}
%%%%%%%%%
In ~\cite{EgOk}, Eguchi and Okamoto
proved the following theorem.
%%%%%%%%%%  Schwartz isomorphism theorem  %%%%%%%%%%%%%
\begin{Thm}[M.\ Eguchi and K.\ Okamoto]
\label{Thm:IsThm}
The Fourier transform $\mathcal{F}$ is a linear topological
isomorphism of $\mathscr{S}(X)$ onto
$\mathscr{I}_{W}(\mathfrak{a}^{*}\times B)$.
\end{Thm}
%%%%%%%%%

Next, we will see some basic properties of Harish-Chandra's
$\boldsymbol{c}$-function and prove a key lemma to apply
pseudodifferential calculi on $\mathfrak{a}$ and $\mathfrak{a}^{*}$
in the later section.

%%%%%%%  Harish-Chandra c-function  %%%%%%%%%%%%%%%%%%%%%%%%%%%%
Let $\Sigma_{0}$ be the indivisible roots of $\Sigma$, that is, whose
elements are $\alpha \in \Sigma$ so that $\alpha / 2 \notin \Sigma$, 
and put $\Sigma_{0}^{+}=\Sigma_{0} \cap \Sigma^{+}$, 
called positive indivisible roots.
Then Harish-Chandra's 
$\boldsymbol{c}$-function can be extended to the
meromorphic function on $\mathfrak{a}^{*}_{\mathbb{C}}$
and is explicitly given by as follows. 
%%%% explicit formula for Harish-Chandra c %%%%%%%%%%%%%
\begin{Thm}[The Gindikin-Karpelevi\v{c} formula]
\label{Thm:HCc}
The $\boldsymbol{c}$-function for the semisimple Lie group $G$
is given by the absolutely convergent integral
$$
 \boldsymbol{c}(\lambda)
=\int_{\Bar{N}}
      e^{-(i\lambda+\rho)(H(\Bar{n}))}
 d\Bar{n},
\quad \operatorname{Re}(i\lambda)
      \in \mathfrak{a}_{+}^{*},
$$
where $\Bar{N}=\theta N$, and 
$d\Bar{n}$ is the Haar measure on
$\Bar{N}$ normalized by
$
\int_{\Bar{N}}
  e^{-2\rho (H(\Bar{n}))}
 d\Bar{n}
=1.
$

Also,
$\boldsymbol{c}(\lambda)$
is explicitly given by the formula
$$
  \boldsymbol{c}(\lambda)
= c_{0}\prod_{\alpha \in \Sigma_{o}^{+}}
  \frac{2^{-\langle i\lambda,\alpha_{0}\rangle}
        \Gamma(\langle i\lambda,\alpha_{0}\rangle)}
       {\Gamma (\frac{1}{2}
        (\frac{1}{2}m_{\alpha}+1+
        \langle i\lambda,\alpha_{0}\rangle))
        \Gamma (\frac{1}{2}
        (\frac{1}{2}m_{\alpha}+m_{2\alpha}
        +\langle i\lambda,\alpha_{0}\rangle))},
$$
where $\alpha_{0}=\alpha /\langle \alpha ,\alpha \rangle$,
the constant $c_{0}$ is given so that 
$\boldsymbol{c}(-i\rho)=1$,
and $\Gamma(z)$ is the gamma function.
\end{Thm}
%%%%%%%%%

%%%%%%%%% Remark for HC c function %%%%%%%%%%%%%%%%%%%%%%%%%%
\begin{Rem}
We have
 \begin{enumerate}[(i)]
 \item Zeros of $\boldsymbol{c}(\lambda)^{-1}$ on $\mathfrak{a}^{*}$
       is precisely the Weyl walls,
       i.e.\ $\bigcup_{\alpha \in \Sigma_{0}^{+}}\ker H_{\alpha}$.
 \item $|\boldsymbol{c}(\lambda)|^{-2}=\boldsymbol{c}(s\lambda)
       \boldsymbol{c}(-s\lambda)$ for all $s\in W,
       \lambda \in \mathfrak{a}^{*}$. 
 \item If $\lambda ,\mu\in\Sigma$ are proportional, i.e.\  
       $\mu=c\lambda$ for some $c\in\mathbb{R}$, then $c=\pm1/2, \pm1$,
       or $\pm2$. $($See 
       ~\cite[Chapter X, Section 3]{HeDLS}.$)$
 \end{enumerate}
\end{Rem}
%%%%%%%%%
In the following, we will use the
pseudodifferential operator theory freely. 
See e.g.\ ~\cite{Ho-III}, ~\cite{Kumano-go}, and ~\cite{Shubin}
for the detail.

A key lemma in this paper is the following.
%%%%% Lemma for HC c-function %%%%%%%%%%%%%%%%%%%%%%%%%%
\begin{Lem}
 \label{Lem:HCcfc}
We have
 \begin{enumerate}
  \item[$(\mathrm{i})$] $\boldsymbol{c}^{-1}(\lambda) \in S^{\dim N/2}$,
               i.e.\ for any $\alpha\in\mathbb{Z}_{\geq 0}^{l}$ there
               exists a positive constant $C_{\alpha} >0$ such that
               $$
               \left|
                   \partial^{\alpha}_{\lambda}
                   \boldsymbol{c}^{-1}(\lambda)
               \right|
               \leq
               C_{\alpha}
               \langle\lambda\rangle^{\dim N /2-|\alpha|}
               $$
               uniformly on $\mathfrak{a}^{*}$.\vspace{-3pt}
  \item[$(\mathrm{ii})$] If $X$ is real rank one, then
               $\boldsymbol{c}^{-1}(\lambda)$ becomes an elliptic symbol
               of order $\dim N/2$.
 \end{enumerate}
\end{Lem}
%%%%%%%%
%%%%%%%%%%% Proof of Lemma 2.7 %%%%%%%%%%%%%%%%%%%%%%%%%%%%%%%%
\begin{proof}
Our proof is almost all same as in
~\cite[Lemma 3.6, p.110]{HeGASS} , 
but we need a little bit more 
precise estimates for $\boldsymbol{c}^{-1}(\lambda)$
to use pseudodifferential calculus.\\
%%%%%%%%%%%%%%%%%%%%%%%% proof of (i) %%%%%%%%%%%%%%%%%%%%%%
$(\mathrm{i})$ 
For the gamma function following formulas are well known:
\begin{alignat*}{2}
 \Gamma(\frac{z}{2})\Gamma(\frac{z+1}{2})
&=2^{1-z}\sqrt{\pi}\,\Gamma(z)
 \quad
&\text{on}\,\, 
 \mathbb{C}\setminus \mathbb{Z}_{\leq 0},\\
 \Gamma(z+1)
&=z\Gamma(z)
 \quad
&\text{on}\,\, 
 \mathbb{C}\setminus \mathbb{Z}_{\leq 0}.
\end{alignat*}
By using these formula, we can write $\boldsymbol{c}^{-1}(\lambda)$
as the following form 
\begin{align}
\boldsymbol{c}^{-1}(\lambda)
&=  c_{0}^{-1}
         \big(\prod_{\alpha\in\Sigma_{0}^{+}}2\sqrt{\pi}\langle
           i\lambda,\alpha_{0}\rangle
          \big)\nonumber\\
&\times \big(\prod_{\alpha\in\Sigma_{0}^{+}}
          s(\tfrac{\langle \lambda,\alpha_{0}\rangle}{2};\tfrac{1}{4}m_{\alpha}
            +\tfrac{1}{2},\tfrac{1}{2})
          s(\tfrac{\langle \lambda,\alpha_{0}\rangle}{2};\tfrac{1}{4}m_{\alpha}
            +\tfrac{1}{2}m_{2\alpha},1)
         \big) \label{eq:EXHCc}
\end{align}
where we put $s(\xi;a,b)=\Gamma(a+i\xi)/\Gamma(b+i\xi)$ for
$\xi \in\mathbb{R}$ and $a,b>0$. So we first investigate $s(\xi;a,b)$.
Using the formula for the gamma function 
$$
\frac{\Gamma'(z)}{\Gamma(z)}
 =\sum_{m=1}^{\infty}(\frac{1}{m}-\frac{1}{m+z})
  -\gamma-\frac{1}{z},
$$
where 
$\gamma 
=\lim_{m\to\infty} 
 \{(1+\tfrac{1}{2}+\cdots+\tfrac{1}{m})-\log m\}
$
is the Euler constant,
we can find that
\begin{align*}
-i\frac{s'(\xi;a,b)}{s(\xi;a,b)} 
&=  \frac{\Gamma'(a+i\xi)}{\Gamma(a+i\xi)}
   -\frac{\Gamma'(b+i\xi)}{\Gamma(b+i\xi)}\\
&=  \sum_{m=0}^{\infty}
    \frac{(a-b)}{(m+a+i\xi)(m+b+i\xi)}.
\end{align*}
So we put
$$
t(\xi;a,b)
=i\sum_{m=0}^{\infty}\frac{(a-b)}{(m+a+i\xi)(m+b+i\xi)},
$$
then we have
\begin{equation}
s'(\xi;a,b)=s(\xi;a,b)t(\xi;a,b).\label{eq:Fos}
\end{equation}
By an elementary calculation, we get $t(\xi;a,b)\in S^{-1}$, so if we
can verify an asymptotic behavior of $s(\xi;a,b)$ as
$|\xi|\to\infty$, then we obtain its of each derivatives of $s(\xi;a,b)$.  
Now we use the formula for the gamma function
$$
\lim_{|z|\to\infty}\frac{\Gamma(z+c)}{\Gamma(z)}e^{-c\log z}=1,
$$
where $|\!\arg z |\leqslant\pi-\delta$, $0<\delta<\pi$, and
$c\in\mathbb{C}$. 
Since we can write 
$$
s(\xi;a,b)
=\frac{\Gamma(a+i\xi)}{\Gamma(i\xi)}
 \frac{\Gamma(i\xi)}{\Gamma(b+i\xi)}
$$
and $|e^{-(a-b)\log (i\xi)}|=|\xi|^{-(a-b)}$, we get
\begin{equation}
\lim_{|\xi|\to \infty}
 |s(\xi;a,b)||\xi|^{-(a-b)}
=1. \label{eq:LIMs}
\end{equation}
Hence by (\ref{eq:EXHCc}), (\ref{eq:Fos}), and $t(\xi;a,b)\in S^{-1}$, we obtain
 $s(\xi;a,b)\in S^{(a-b)}$.
Since $\langle \lambda,\alpha_{0}\rangle$ is 
a polynomial of $\lambda\in\mathfrak{a}^{*}$ 
of order one, 
by using expression (\ref{eq:EXHCc})
we have $\boldsymbol{c}^{-1}(\lambda)\in S^{m}$.
By a simple cluculation, we have 
$m={\dim N}/2$. 

%%%%%%%%%%%%%%%%%%%%%%%%% proof of (ii) %%%%%%%%%%%%%%%%%%%%%%%%%%%%
$(\mathrm{ii})$ If $X$ is real rank one, then 
$|\langle i\lambda,\alpha_{0}\rangle|=|\lambda|$ for all
$\lambda$, $\alpha\in\mathfrak{a}^{*}$, 
so ellipticity follows from (\ref{eq:EXHCc}) and (\ref{eq:LIMs}).  
 \qedhere
\end{proof}

Finally, we shall introduce the Radon transform and see 
its microlocal properties. 

For $f\in\mathscr{D}(X)$,
the Radon transform is defined by
$$
 \mathcal{R}f(\xi)
=\int_{x\in\xi}
      f(x)
 dm_{\xi}(x)
 \quad \text{for}\,\,\xi\in\Xi , 
$$
where $dm_{\xi}$ denotes the induced measure on $\xi$,
or equivalently under the identification
$\mathfrak{a}\times B\ni (H,b)
 \leftrightarrow \xi(H,b) \in\Xi$,
we can also write
$$
\mathcal{R}f(H,kM)
=\int_{N}f(k\exp(H)n\cdot o)
  dn.
$$
For $\varphi\in L_{\mathrm{loc}}^{1}(\Xi,d\xi)$, the dual Radon
transform is defined by 
$$
\mathcal{R}^{*}\varphi(x)
=\int_{\xi\ni x}
     \varphi(\xi)
 dm_{x}(\xi)
\quad \text{for}\,\, x\in X ,
$$
where $dm_{x}$ denotes the induced measure on
$\{ \xi\in\Xi ;\xi\ni x\}$.
Then $\mathcal{R}^{*}$ is the formal adjoint of $\mathcal{R}$
in the following sense:
$$
 \int_{X}
  f(x)\mathcal{R}^{*}\varphi(x)
  dx
=\int_{\Xi}
  \mathcal{R}f(\xi)\varphi(\xi)
  d\xi,
\quad f\in\mathscr{D}(X),\varphi\in\mathscr{E}(\Xi).
$$

For the Radon transform, the following inversion formula holds.
(For the detail and the proof, see e.g.\ ~\cite{HeGASS}.)
%%%%%%%%%% Radon inversion formula %%%%%%%%%%%%%%%%%
\begin{Thm}
 \label{Thm:Radoninv}
For any $u\in\mathscr{S}(X)$, the Radon transform is 
invertible by
$$
 u
=w^{-1}\mathcal{R}^{*}
       \overline{\Lambda}
       \Lambda
       \mathcal{R}
\quad\text{on}\quad X.
$$
Where $\Lambda$ is a pseudodifferential operator on 
$\Xi$ defined by 
$$
 \Lambda
=e^{-\rho (H)}
 \boldsymbol{c}^{-1}(D_{H})
 e^{\rho(H)}.
$$
\end{Thm} 
%%%%%%%%%

Next, we will see microlocal properties of the
Radon transform and establish a key theorem.

In ~\cite{GuilSter} Guillemin and Sternberg introduced 
microlocal techniques to the study of generalized Radon 
transform between two manifolds $X$, $Y$ by noting 
that $\delta_{Z}$, where $Z\subset X\times Y$ is a submanifold
called incidence relation, is an example of a Fourier integral
distribution and proceeding to study the microlocal analogue 
of the double fibration. 
After that, 
Beylkin ~\cite{Beylkin}, 
Greenleaf and Uhlmann ~\cite{GreUhl}, 
Quinto ~\cite{QuintoRadon},
Gonzalez and Quinto ~\cite{GonzalezQuinto} 
developed the study of microlocal properties, 
and applied to the more general support theorems,
partial differential equations, etc. 

Let $Z \subset\Xi\times X$ be the incidence relation, 
that is, 
$
Z=\{(\xi ,x)\in\Xi\times X;x\in\xi\},
$
and we consider the following double fibration diagram.

%%%%% Figure of diagram 1 %%%%%%%%%%%%%%%%%%%%%%%%%%%%%%%
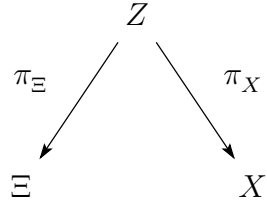
\begin{figure}[h]
\begin{center}
\unitlength 0.1in
\begin{picture}( 16.8000,  9.0000)( 52.2000,-46.6500)
% STR 2 0 3 0
% 3 6400 3750 6400 3850 5 0
% $Z$
\put(64.0000,-38.5000){\makebox(0,0){$Z$}}%
% STR 2 0 3 0
% 3 5800 4650 5800 4750 5 0
% $\Xi$
\put(58.0000,-47.5000){\makebox(0,0){$\Xi$}}%
% STR 2 0 3 0
% 3 7000 4650 7000 4750 5 0
% $X$
\put(70.0000,-47.5000){\makebox(0,0){$X$}}%
% STR 2 0 3 0
% 3 6950 4100 6950 4200 5 0
% $\pi^{}_{X}$
\put(69.5000,-42.0000){\makebox(0,0){$\pi^{}_{X}$}}%
% STR 2 0 3 0
% 3 5850 4100 5850 4200 5 0
% $\pi^{}_{\Xi}$
\put(58.5000,-42.0000){\makebox(0,0){$\pi^{}_{\Xi}$}}%
% VECTOR 2 0 3 0
% 4 6300 4000 5900 4600 6500 4000 6900 4600
% 
\special{pn 8}%
\special{pa 6300 4000}%
\special{pa 5900 4600}%
\special{fp}%
\special{sh 1}%
\special{pa 5900 4600}%
\special{pa 5954 4556}%
\special{pa 5930 4556}%
\special{pa 5920 4534}%
\special{pa 5900 4600}%
\special{fp}%
\special{pa 6500 4000}%
\special{pa 6900 4600}%
\special{fp}%
\special{sh 1}%
\special{pa 6900 4600}%
\special{pa 6880 4534}%
\special{pa 6870 4556}%
\special{pa 6846 4556}%
\special{pa 6900 4600}%
\special{fp}%
\end{picture}%
\end{center}
 \end{figure}
%%%%%%%%%%%%%%%%%%%%%%%%%%%%%%%%%%%%%%%%%%%%%%%%%%%%
\vspace{10pt}
Then the Schwartz kernel of $\mathcal{R}$ is $\delta_{Z}$
, that is, the delta function supported on $Z$,
and $\delta_{Z}$ is a Fourier integral distribution,
so from H\"{o}rmander's theory Radon transform
$\mathcal{R}$ is a Fourier integral operator of
order $-\dim N/2$ associated with the canonical
relation $\Gamma =N^{*}Z\setminus 0$, the twisted
conormal bundle of $Z$.
That is,
\begin{align*}
\Gamma
 &=\left\{ (\xi,\eta_{\xi};x,\omega_{x})
            \in T^{*}(\Xi\times X)\setminus 0\,
            ;(\xi,x)\in Z,
             (\eta_{\xi},-\omega_{x})
             \bot T_{(\xi,x)}Z\right\}\\
 &=\left\{ \left( H,b,x;
                  d_{x}(\lambda(A(x,b)))
                 +d_{H}(\lambda(H))
                 -d_{b}(\lambda(A(x,b)))
           \right) 
           ;\lambda\in\mathfrak{a}^{*}\setminus \{0\},
            H=A(x,b)
   \right\}.
\end{align*}
Similarly, $\mathcal{R}^{*}$ is a Fourier integral operator
associated with the canonical relation
$\Gamma^{t} \subset T^{*} X\times T^{*}\Xi$,
which is simply $\Gamma$ with $(x,\omega_{x})$ and
$(\xi,\eta_{\xi})$ interchanged.

\newpage

Now consider the microlocal diagram

%%%%%%%%%%% Figure of diagram 2 %%%%%%%%%%%%%%%%%%%%%%%%%%%%%
\begin{figure}[h]
\begin{center}
\unitlength 0.1in
\begin{picture}( 19.0500,  9.0000)( 49.9500,-46.6500)
% STR 2 0 3 0
% 3 6400 3750 6400 3850 5 0
% $\Gamma$
\put(64.0000,-38.5000){\makebox(0,0){$\Gamma$}}%
% STR 2 0 3 0
% 3 5800 4650 5800 4750 5 0
% $T^{*}\Xi$
\put(58.0000,-47.5000){\makebox(0,0){$T^{*}\Xi$}}%
% STR 2 0 3 0
% 3 7000 4650 7000 4750 5 0
% $T^{*}X$
\put(70.0000,-47.5000){\makebox(0,0){$T^{*}X$}}%
% STR 2 0 3 0
% 3 6950 4100 6950 4200 5 0
% $\pi^{}_{T^{*}X}$
\put(69.5000,-42.0000){\makebox(0,0){$\pi^{}_{T^{*}X}$}}%
% STR 2 0 3 0
% 3 5850 4100 5850 4200 5 0
% $\pi^{}_{T^{*}\Xi}$
\put(58.5000,-42.0000){\makebox(0,0){$\pi^{}_{T^{*}\Xi}$}}%
% VECTOR 2 0 3 0
% 4 6300 4000 5900 4600 6500 4000 6900 4600
% 
\special{pn 8}%
\special{pa 6300 4000}%
\special{pa 5900 4600}%
\special{fp}%
\special{sh 1}%
\special{pa 5900 4600}%
\special{pa 5954 4556}%
\special{pa 5930 4556}%
\special{pa 5920 4534}%
\special{pa 5900 4600}%
\special{fp}%
\special{pa 6500 4000}%
\special{pa 6900 4600}%
\special{fp}%
\special{sh 1}%
\special{pa 6900 4600}%
\special{pa 6880 4534}%
\special{pa 6870 4556}%
\special{pa 6846 4556}%
\special{pa 6900 4600}%
\special{fp}%
\end{picture}% 
\end{center}
\end{figure}
%%%%%%%%%%%%%%%%%%%%%%%%%%%%%%%%%%%%%%%%%%%%%%%%%%%%%%%%%%%
\noindent
where $\pi^{}_{T^{*}\Xi}$ and $\pi^{}_{T^{*}X}$ again
denote the natural projections.
If $X$ is real rank one, 
then $\pi_{\Xi}$ becomes an
injective immersion (cf.~\cite{QuintoRadon}),
i.e.\ satisfies the Bolker assumption,
so $\Gamma\subset T^{*}\Xi\times T^{*}X$
becomes a local canonical graph.
Then by a general Fourier integral
operator theory and note for Sobolev wave front
set in ~\cite{QuintoXray}, we have the following theorem.
%%%%%%%%%%% Continuity of the Radon transform %%%%%%%%%%%%%%%%%%
\begin{Thm}
 \label{Thm:SoWF}
Let $X$ be a real rank one symmetric space of noncompact type.
Then we have the followings:
 \begin{itemize}
  \item[$\mathrm{(i)}$] For any $s\in\mathbb{R}$, 
               we have the following continuous maps:
               \begin{align*}
                \mathcal{R}
                :H_{\mathrm{comp}}^{s}(X)
                 \rightarrow
                 H_{\mathrm{loc}}^{s+\dim N/2}(\Xi),\\
                \mathcal{R}^{*}
                :H_{\mathrm{comp}}^{s}(\Xi)
                 \rightarrow
                 H_{\mathrm{loc}}^{s+\dim N/2}(X).
                \end{align*}
  \item[$\mathrm{(ii)}$] For any $u\in\mathscr{E}'(X)$,
               $\psi\in\mathscr{E}'(\Xi)$,
               and $s\in\mathbb{R}$ we have
               \begin{align*}
                \operatorname{WF}^{s+\dim N/2}
                (\mathcal{R}u)
                &\subset
                \Gamma\circ
                \operatorname{WF}^{s}(u),\\
                \operatorname{WF}^{s+\dim N/2}
                (\mathcal{R}^{*}\psi)
                &\subset
                \Gamma^{t}\circ
                \operatorname{WF}^{s}(\psi),
                \end{align*}
               where for $C\subset T^{*}\Xi\times T^{*}X$,
               $A\subset T^{*}X$, $B\subset T^{*}\Xi$,
               put
               \begin{align*}
                C\circ
                A&=\{ (\xi;\eta_{\xi})\in T^{*}\Xi
                     ;(\xi;\eta_{\xi},x;\omega_{x})
                      \in C \,\,\text{for\,some}\,
                      (x;\omega_{x})\in A\}, \\
                C^{t}\circ
                B&=\{ (x;\omega_{x})\in T^{*}X
                     ;(x;\omega_{x},\xi;\eta_{\xi})
                      \in C^{t} \,\,\text{for\,some}\,
                      (\xi;\eta_{\xi})\in B \}.
               \end{align*}
 \end{itemize}
\end{Thm}
%%%%%%%%%%%%%%%%%%%%%%%%%%%%%%%%%%%%%%%%%%%%%%%%%%%%%%%%%%%%%%%%%%%%%%
%  Section3 : Basic properties of weighted L2space                   %
%%%%%%%%%%%%%%%%%%%%%%%%%%%%%%%%%%%%%%%%%%%%%%%%%%%%%%%%%%%%%%%%%%%%%%
\section{Basic properties of weighted $L^{2}$-space\label{WL2}}
This section is devoted to proving local regularities of
weighted $L^{2}$-space.

We shall study the properties of weighted $L^{2}$-space 
$L^{2,\delta}(X)$, by using  pseudodifferential calculi 
on $\mathfrak{a}$, $\mathfrak{a}^{*}$ 
and the Radon inversion formula.
%%%% Proposition for well-definedness of w-L2 space %%%%%%%%%%
\begin{Prop}
 \label{Prop:WDwL2}
For any $\delta\in\mathbb{R}$, there exists a continuous seminorm
$\| \cdot\|$ on $\mathscr{S}(X)$ such that for all $u\in\mathscr{S}(X)$
$$
 \|u \|_{L^{2,\delta}(X)}
=\| \langle D_{\lambda}\rangle^{\delta}
 (\mathcal{F}u(\lambda,b)\boldsymbol{c}(\lambda)^{-1})
  \|_{L^{2}(\mathfrak{a}^{*}\times B, w^{-1}d\lambda db)}
 \leqslant \|u\| .
$$
In particular, $\|u\|_{ L^{2,\delta} (X) } < \infty$ for all 
$u\in \mathscr{D} (X)$, hence $L^{2,\delta}(X)$ is the 
well-defined Hilbert space.
\end{Prop}
%%%%%%%%%%
%%%%%%%%%%% Proof of well definedness %%%%%%%%%%%%%%%%%%%%%%%
\begin{proof}
By using Lemma \ref{Lem:HCcfc} and Theorem \ref{Thm:Radoninv}, 
the assertion is obvious.
\end{proof}
%%%%%%%%%%%

The basic properties of the weighted $L^{2}$-space
are the following.
%%%%% Proposition for basic properties of W-L2 space %%%%%%%%%%
\begin{Prop}
 \label{Prop:BPwL2}
 We have
\begin{itemize}
\item[$(\mathrm{i})$] 
 $L^{2,0}(X)=L^{2}(X).$\vspace{-1mm}
\item[$(\mathrm{ii})$] 
 $L^{2,\delta}(X)\hookrightarrow L^{2,\delta'}(X)$ 
 \quad $(\delta\geqslant\delta').$\vspace{-1mm}
\item[$(\mathrm{iii})$]
 For $g\in G$, put $\tau_{g}u(x)=u(g\cdot x)$.
 Then for any $\delta\in\mathbb{R}$ 
 we have the following linear continuous map:
 $$
 L^{2,\delta}(X)
 \ni
 u
 \mapsto
 \tau_{g}u
 \in 
 L^{2,\delta}(X).
 $$
 Hence $L^{2,\delta}(X)$ is a 
 $G$-invariant Hilbert space.
\item[$(\mathrm{iv})$]
 For any $\delta\in\mathbb{R}$, the following
 inclusion maps are all continuous dense embedding:
 $$
 \mathscr{D}(X)
 \hookrightarrow \mathscr{S}(X)
 \hookrightarrow L^{2,\delta}(X).
 \vspace{-2mm}
 $$
\item[$(\mathrm{v})$]
 By the natural coupling we have
 $$
 L^{2,\delta}(X) \hookrightarrow
 \mathscr{S}'(X).
 $$
\end{itemize}
\end{Prop}
%%%%%%%%%%
%%%%%%%%%%%% Proof of basic property %%%%%%%%%%%%%%%%%%%%%%%%%%%
\begin{proof} 
%%%%%%%%%%%%%%%%%%%%%%%%%% proof of (i) %%%%%%%%%%%%%%%%%%%%%
$(\mathrm{i})$ 
It follows from the Plancherel theorem immediately.\\
%%%%%%%%%%%%%%%%%%%%%%%%%% proof of (ii) %%%%%%%%%%%%%%%%%%%%
$(\mathrm{ii})$ 
Since for any $s>0$, 
$\langle D_{\lambda}\rangle^{-s}:
 L^{2}(\mathfrak{a}^{*})\rightarrow
 L^{2}(\mathfrak{a}^{*})$
is continuous, it's obvious.\\
%%%%%%%%%%%%%%%%%%%%%%%%%% proof of (iii) %%%%%%%%%%%%%%%%%%%%%%%%%
$(\mathrm{iii})$ For any $g\in G$ and $u\in\mathscr{D}(X)$, 
we have 
$$
\mathcal{F}(\tau_{g}u)(\lambda, b)
=e^{(i\lambda-\rho)(A(g\cdot o, g\cdot b))}
 \mathcal{F}u(\lambda, g\cdot b).
$$
From the definition of the weighted $L^{2}$-norm, 
we see that
\begin{align*}
& \|
   \tau_{g}u
  \|_{L^{2,\delta}(X)}\\
=&\|
   \langle D_{\lambda}\rangle^{\delta}
   \left( 
    \mathcal{F}(\tau_{g}u)
               (\lambda, b)
    \boldsymbol{c}^{-1}(\lambda)
   \right)
  \|_{L^{2}
      (\mathfrak{a}^{*}
      \times B,
      w^{-1}dHdb)}\\
=&\|
   \left(
     \langle D_{\lambda}\rangle^{\delta}
     e^{(i\lambda-\rho)(A(g\cdot o,g\cdot b))}
     \langle D_{\lambda}\rangle^{-\delta}
    \right)
   \langle D_{\lambda}\rangle^{\delta}
   \left( 
     \mathcal{F}u
               (\lambda, b)
     \boldsymbol{c}^{-1}(\lambda)
    \right)
  \|_{L^{2}
      (\mathfrak{a}^{*}
      \times B,
      w^{-1}dHdb)}.
\end{align*}
By using the basic theory of pseudodifferential
operators, we can find 
$$
\sup_{b\in B}
\|
  \langle D_{\lambda}\rangle^{\delta}
  e^{(i\lambda-\rho)(A(g\cdot o, g\cdot b))}
  \langle D_{\lambda}\rangle^{-\delta}
\|_{\mathcal{L}
    (L^{2}(\mathfrak{a}^{*}))}
< \infty ,
$$
where 
$
\|\cdot\|_{\mathcal{L}
(L^{2}(\mathfrak{a}^{*})}
$
is a operator norm on $L^{2}(\mathfrak{a}^{*})$.
Thus we obtain
$$
\|
   \tau_{g}u
 \|_{L^{2,\delta}(X)}
\leq
\sup_{b\in B}
 \|
   \langle D_{\lambda}\rangle^{\delta}
   e^{(i\lambda-\rho)(A(g\cdot o, g\cdot b))}
   \langle D_{\lambda}\rangle^{-\delta}
 \|_{\mathcal{L}
    (L^{2}(\mathfrak{a}^{*}))}
 \|u\|_{L^{2,\delta}(X)}.
$$
Since $\mathscr{D}(X)$ is dense in $L^{2,\delta}(X)$, 
we get the linear continuous map.\\
%%%%%%%%%%%%%%%%%%%%%%%%%% proof of (iv) %%%%%%%%%%%%%%%%%%%%%%%%%%
$(\mathrm{iv})$
By the definition of $L^{2,\delta}(X)$ and
from Proposition \ref{Prop:WDwL2}, it holds.\\ 
%%%%%%%%%%%%%%%%%%%%%%%%% proof of (v) %%%%%%%%%%%%%%%%%%%%%
$(\mathrm{v})$
We can easily see that the follwing Hermitian form
on $\mathscr{D}(X)\times\mathscr{D}(X)$:
$$
 (u,v)
=w^{-1}
 \iint_{\mathfrak{a}^{*}\times B}
   \langle D_{\lambda}\rangle^{\delta}
   \big( 
         \mathcal{F}u(\lambda ,b)
         \boldsymbol{c}^{-1}(\lambda) 
   \big) 
  \overline
  { \langle D_{\lambda}\rangle^{\delta}
    \big( 
          \mathcal{F}v(\lambda ,b)
          \boldsymbol{c}^{-1}(\lambda) 
    \big) }
 d\lambda db,
$$
can be uniquely continuously extended to its on
$L^{2,\delta}(X)\times\mathscr{S}(X)$.
\end{proof}
%%%%%%%%%%%
Now, we will examine the local regularity of the
weighted $L^{2}$-space $L^{2,\delta}(X)$.

%%%%% Proposition for regularity of w-L2 space %%%%%%%%
\begin{Prop}
 \label{Prop:REwL2}
Suppose $X$ is real rank one, then we have 
following continuous embeddings:
 \begin{equation}
  L^{2}_{\mathrm{comp}}(X)
  \hookrightarrow
  L^{2,\delta}(X)
  \hookrightarrow
  L^{2}_{\mathrm{loc}}(X).
 \end{equation} 
\end{Prop}
%%%%%%%%%%
%%%%% Proof of local regularity %%%%%%%%%%%%%%%%%%%%%%%%%%%%%%%%%%%%%%%%%%%%%%
\begin{proof}
%%%%% first inclusion %%%%%%%%%%%%%%%%%%%%%
First, we will prove 
$
L^{2}_{\mathrm{comp}}(X)
\hookrightarrow
L^{2,\delta}(X).
$
It is enough to show for any $R>0$ and 
$k \in\mathbb{N}$ that there exists
a constant $C_{k,R}>0$ such that for all
$u\in\mathscr{D}(X)$ with 
$\operatorname{supp} u\subset \overline{B(o,R)}$
the following estimate holds
$$
\|u\|_{L^{2,2k}(X)}\leq C_{k,R}\|u\|_{L^{2}(X)}.
$$
Take $\chi_{1}\in\mathscr{D}(\mathfrak{a})$ with
$0\leq\chi_{1}\leq 1$ such that $\chi_{1} =1$ near 
$\{H\in\mathfrak{a};|H|\leq R\}$.
Then for any $j\in\mathbb{N}$, we have 
\begin{align*}
&\frac{\partial^{j}}{\partial \lambda^{j}}
  (\mathcal{F}u(\lambda ,b)\boldsymbol{c}^{-1}(\lambda))\\
&=\sum_{j_{1}+j_{2}=j}
  \frac{j!}{j_{1}!j_{2}!}
  \Big( \frac{\partial^{j_{2}}}{\partial \lambda^{j_{2}}}
         \boldsymbol{c}^{-1}(\lambda ) 
   \Big)
  \int_{\mathfrak{a}}e^{-i\lambda (H)}(-iH)^{j_{1}}
   e^{\rho (H)}\mathcal{R} u(H,b)dH;\\
\intertext{ By putting
            $p_{j_{2}}(\lambda)
             =\partial^{j_{2}}_{\lambda}
            \boldsymbol{c}^{-1}(\lambda )$, we can rewrite}
&=\int_{\mathfrak{a}}e^{-i\lambda (H)}
        \sum_{j_{1}+j_{2}=j}\frac{j!}{j_{1}!j_{2}!}
        \Bar{p}_{j_{2}}(D_{H})^{*}
        (-iH)^{j_{1}}\chi_{1}(H)
         \left( e^{\rho (H)}\mathcal{R}u(H,b)
         \right) dH\\
&=\int_{\mathfrak{a}}e^{-i\lambda (H)}q_{j}(H,D_{H})
        \left(e^{\rho (H)}\mathcal{R}u(H,b)\right)dH,
\end{align*}
where we set 
$$
 q_{j}(H,D_{H})
=\sum_{j_{1}+j_{2}=j}
 \frac{j!}{j_{1}!j_{2}!}
 \Bar{p}_{j_{2}}(D_{H})^{*}
 (-iH)^{j_{1}}\chi_{1}(H) .
$$
Since $q_{j}(H,D_{H})\in OpS^{\dim N/2}$,
there exists a constant 
$C_{R,j}>0$ such that
\begin{align*}
\left\|
        \frac{\partial^{j}}{\partial \lambda^{j}}
        (\mathcal{F}u(\lambda ,b)\boldsymbol{c}^{-1}(\lambda))
\right\|_{L^{2}(\mathfrak{a}^{*}\times B,w^{-1}d\lambda db)}
&\leq C_{R,j}
  \| e^{\rho (H)\mathcal{R}u(H,b)}
    \|_{L^{2}(B;H^{\dim N/2}(\mathfrak{a}))}\\
&=C_{R,j}
  \| \langle\lambda\rangle^{\dim N/2}\mathcal{F}u(\lambda ,b)
    \|_{L^{2}(\mathfrak{a}^{*}\times B, d\lambda db)}.
\end{align*}
By the assumption, from Lemma \ref{Lem:HCcfc} (ii),
$\boldsymbol{c}^{-1}(\lambda)$ is the elliptic symbol of 
order $\dim N/2$, so there exist $R'>0$ and $C'>0$ such that 
for all $|\lambda |\geq R'$ 
$$
|\boldsymbol{c}^{-1}(\lambda)|
 \geq {C'}^{-1}
      \langle \lambda \rangle^{\dim N/2}.
$$
Take $\chi_{2}\in\mathscr{D}(\mathfrak{a})$
with $0\leq\chi_{2}\leq 1$, 
$\operatorname{supp} \chi_{2}
 \subset 
 \{H\in\mathfrak{a};|H|\leq R'+1\}$ 
and 
$\chi_{2}=1$ near $\{H\in\mathfrak{a};|H|\leq R'+1\}$, 
then we have 
\begin{align*}
&      \left| \langle\lambda\rangle^{\dim N/2}
        \mathcal{F}u(\lambda ,b)\right|\\
&=     \langle\lambda\rangle^{\dim N/2}
       \chi_{2}(\lambda)\left|\mathcal{F}u(\lambda,b)\right|
      +\langle\lambda\rangle^{\dim N/2}
       (1-\chi_{2}(\lambda))
       \left|\mathcal{F}u(\lambda,b)\right|\\
&\leq \langle R'\rangle^{\dim N/2}
       \chi_{2}(\lambda)\left|\mathcal{F}u(\lambda,b)\right|
      +C'|\boldsymbol{c}^{-1}(\lambda)|
       (1-\chi_{2}(\lambda))
       \left|\mathcal{F}u(\lambda,b)\right|\\
&\leq \langle R'\rangle^{\dim N/2}
       \chi_{2}(\lambda)\left|\mathcal{F}u(\lambda,b)\right|
      +C'|\mathcal{F}u(\lambda,b)
         \boldsymbol{c}^{-1}(\lambda)|.
\end{align*}
For the first term in above inequality,
by using the H\"{o}lder inequality  we find
\begin{align*}
  \left|\mathcal{F}u(\lambda,b)\right|
&=    \left|\int_{X}e^{(-i\lambda+\rho)(A(x,b))}
            u(x)dx\right| \\
&\leq \int_{\overline{B(o,R)}}e^{\rho (A(x,b))}
           |u(x)|dx \\
&\leq e^{|\rho |R}
      \operatorname{Vol}(B(o,R))^{\frac12}\|u\|_{L^{2}(X)}.
\end{align*}
By using the Plancherel formula, there exists a constant
$C_{R,R',j}>0$ such that
$$
\left\| \frac{\partial^{j}}{\partial \lambda^{j}}
        (\mathcal{F}(\lambda ,b)\boldsymbol{c}^{-1}(\lambda))
 \right\|_{L^{2}(\mathfrak{a}^{*}\times B,w^{-1}\lambda db)}
\leq C_{R,R',j} \|u\|_{L^{2}(X)}.
$$
Since 
$\langle D_{\lambda}\rangle^{2k}
 =(1-\partial^{2}_{\lambda})^{k}$,
we obtain desired estimates.

%%%%%%%%% second inclusion  %%%%%%%%%%%%%%%%%%%%
\vspace{4pt}
Next, we will prove
$L^{2,\delta}(X)
 \hookrightarrow
 L^{2}_{\mathrm{comp}}(X).$
For any $u\in\mathscr{D}(X)$ and $\chi_{1}\in\mathscr{D}(X)$, 
by using the Radon inversion formula we have
$$
\chi_{1}u=w^{-1}\chi_{1}
          \mathcal{R}^{*}\Bar\Lambda
          \Lambda\mathcal{R}u.
$$
Since $\operatorname{supp}\chi_{1}$ is compact, 
there exists a sufficiently large $R>0$ such that
$\operatorname{supp}\chi_{1}\subset\overline{B(o,R)}$.
We can take $\chi_{2}\in\mathscr{D}(\mathfrak{a})$ 
with $0\leq \chi_{2}\leq 1$,
$\operatorname{supp}\chi_{2} \subset
 \{H \in \mathfrak{a} ;|H|\leq R+1\}$
and $\chi_{2}=1$ near $\{H\in\mathfrak{a};|H|\leq R\}$. 
Then we have
$$
\mathcal{R}^{*}\chi_{2}\Bar\Lambda\Lambda\mathcal{R}u
= \mathcal{R}^{*}\Bar\Lambda \Lambda\mathcal{R}u
$$
on $\overline{B(o,R)}$,
hence we obtain
\begin{equation*}
 \chi_{1}u
=w^{-1}\chi_{1}
 \mathcal{R}^{*}\chi_{2}\Bar\Lambda
 \Lambda\mathcal{R}u 
 \quad\text{in}\,\,
 \mathscr{D}(X).
\end{equation*}
By a density argument, 
we can uniquely extend the linear map 
$$
\mathscr{D}(X)
\ni u
\mapsto
\langle H\rangle^{\delta}
\Lambda
\mathcal{R}u
\,\,\in
L^{2}(\Xi ,w^{-1}d\xi) 
$$
to the linear isometry operator
$$
L^{2,\delta}(X)
\ni u
\mapsto
\varphi_{u}
\,\,\in
L^{2}(\Xi ,w^{-1}d\xi) 
$$
which satisfying
$$
  \chi_{1}u 
= w^{-1}
  \chi_{1}\mathcal{R}^{*}\chi_{2}
  \Bar\Lambda
  \langle H\rangle^{-\delta}
  \varphi_{u}
  \quad\text{in}\,\,
  \mathscr{E}'(X)
$$
for any $u\in L^{2,\delta}(X)$.

Now take $\chi_{3}\in\mathscr{D}(\mathfrak{a})$ such that 
$0\leq\chi_{3}\leq1$ and $\chi_{3}=1$ near 
$\operatorname{supp}\chi_{2}$,
then we have  
$\operatorname{supp}\chi_{2}\cap
 \operatorname{supp} (1-\chi_{3})
=\emptyset$. 
Then we can decompose $\chi_{1}u\in\mathscr{E}'(X)$ as
\begin{align*}
\chi_{1}u &= w^{-1}\chi_{1}\mathcal{R}^{*}\chi_{2}\Bar\Lambda
                   \chi_{3}\langle H\rangle^{-\delta}
                   \varphi_{u}
            +w^{-1}\chi_{1}\mathcal{R}^{*}\chi_{2}\Bar\Lambda
                   (1-\chi_{3})\langle H\rangle^{-\delta}
                   \varphi_{u}\\
          &=u_{1}+u_{2}.
\end{align*}
To estimate the first term, set 
$$
q(H,D_H)
=\chi_{2}
 e^{-\rho (H)}
 \overline{\boldsymbol{c}}^{-1}(D_{H})
 \chi_{3}\langle H\rangle^{-\delta}
 \in OpS^{\operatorname{dim} N/2}.
$$
Then we have 
$$
u_{1}
 =w^{-1}\chi_{1}\mathcal{R}^{*}p(H,D_{H})
  \left( e^{\rho (H)}\varphi_{u} \right).
$$
Since all the following maps are continuous
$$
\begin{array}{rclcl}
p(H,D_{H})
 \!\!\!&:&\!\!\! L^{2}           (\mathfrak{a}\times B,dHdb)
 \!\!\!&\rightarrow &\!\!\!      L^{2}(B;H_{\mathrm{comp}}^{-\dim N/2}
                                 (\mathfrak{a})),\vspace{1mm}\\
\iota
 \!\!\!&:&\!\!\!                 L^{2}(B;H_{\mathrm{comp}}^{-\dim N/2}(\mathfrak{a}))
 \!\!\!&\hookrightarrow &\!\!\!  H_{\mathrm{comp}}^{-\dim N/2}
                                 (\mathfrak{a}\times B),\vspace{1mm}\\
\mathcal{R}^{*}
 \!\!\!&:&\!\!\!                 H_{\mathrm{comp}}^{-\dim N/2}(\mathfrak{a}\times B)
 \!\!\!&\rightarrow &\!\!\!      L_{\mathrm{loc}}^{2}(X),
 \end{array}
$$
there exists a constant $C_{1}>0$ such that
\begin{align}
\|u_{1}\|_{L^{2}}
 &\leq
  C_{1}
  \| e^{\rho (H)}\varphi_{u}(H,b)
   \|_{L^{2}(\mathfrak{a}\times B,w^{-1}dHdb)} \nonumber\\
 &\leq
  C_{1}
  \|\varphi_{u}
   \|_{L^{2}(\Xi , w^{-1}d\xi)}. \label{eq:ESu1}
\end{align}
For the second term, set
$$
r(H,D_{H})
=\chi_{2}
 e^{-\rho (H)}
 \overline{\boldsymbol{c}}^{-1}(D_{H})
 (1-\chi_{3})\langle H\rangle^{-\delta}
 \in OpS^{-\infty},
$$
then we can write 
$$
u_{2}
 =w^{-1}\chi_{-1}\mathcal{R}^{*}r(H,D_{H})
  \left( \langle H\rangle^{-\delta }
       e^{\rho (H)}\varphi_{u} \right).
$$
Since $r(H,D_{H})\in OpS^{-\infty}$, there exists the kernel 
$K(H,H')\in \mathscr{E}(\mathfrak{a}\times\mathfrak{a})$ 
which satisfies for any $j_{1},j_{2},j_{3}\in\mathbb{N}$ 
$$
\sup_{(H,H')\in\mathfrak{a}\times\mathfrak{a}}
 \left| \langle H'\rangle^{-j_{1}}
        \partial_{H}^{j_{2}}
        \partial_{H'}^{j_{3}}
        K(H,H') \right|
 <\infty ,
$$
and for any $\varphi \in \mathscr{S}'(\mathfrak{a})$
$$
r(H,D_{H})\varphi (H)
=\int_{\mathfrak{a}  }K(H,H-H')\varphi (H')dH'.
$$
Then the distribution kernel of 
$r(H,D_{H})\langle H\rangle^{-\delta}
 =\chi_{3} r(H,D_{H})\langle H\rangle^{-\delta}$ 
is
$$
\chi_{3}(H)K(H,H-H')\langle H'\rangle^{-\delta}
\in L^{2}(\mathfrak{a}\times\mathfrak{a},dHdH'),
$$
i.e.\ the Hilbert-Schmidt kernel, so we have the 
following continuous map
$$
r(H,D_{H})\langle H\rangle^{-\delta}:
L^{2}(\mathfrak{a}\times B,dHdb)
 \rightarrow L^{2}_{\mathrm{comp}}(\mathfrak{a}\times B,dHdb),
$$
as well as
$$
\mathcal{R}^{*}: L_{\mathrm{comp}}^{2}(\Xi)
 \rightarrow     H_{\mathrm{loc}}^{\dim N/2}(X)
 \hookrightarrow L_{\mathrm{loc}}^{2}(X).
$$
Hence there exists a constant $C_{2}>0$ such that 
\begin{align}
\|u_{2}\|_{L^{2}}
 &\leq  C_{2}\| e^{\rho (H)}\varphi_{u}(H,b)
              \|_{L^{2}(\mathfrak{a}\times B,w^{-1}dHdb)} \nonumber\\
 &\leq  C_{2}\|\varphi_{u}\|_{L^{2}(\Xi , w^{-1}d\xi)} \label{eq:ESu2}
\end{align}
Thus from (\ref{eq:ESu1}), (\ref{eq:ESu2}), we obtain 
$\chi_{1}u\in L^{2}(X)$ and there exists  a constant $C_{R}>0$ 
such that 
$$
\|\chi_{1}u\|_{L^{2}(X)}
 \leq C_{R}\|\varphi_{u}\|_{L^{2}(\Xi ,w^{-1}d\xi)}
 =    C_{R}\|u\|_{L^{2,\delta}(X)}.
$$ 
So we complete the proof.
\end{proof}
%%%%%%%%%%%
%%%%%%%%%%%%%%%%%%%%%%%%%%%%%%%%%%%%%%%%%%%%%%%%%%%%% 
% Section4: Time-global smoothing effects           %
%%%%%%%%%%%%%%%%%%%%%%%%%%%%%%%%%%%%%%%%%%%%%%%%%%%%%
\section{Time-global smoothing effects\label{TGSE}}
In this section , we prove Theorem \ref{Thm:SmEf}, that is, 
a smoothing effect of radially symmetric constant 
coefficients pseudodifferential equations on 
real rank one symmetric spaces of noncompact type. 
The main idea of the proof is to reduce the argument to 
the Euclidean case by introducing some isometry 
operator from $L^{2}$ space on a Riemannian symmetric space 
to its on the horocycle space. 
We will see later this isometry transform a solution on the 
Riemannian symmetric space  to the solution of the ``same'' equation
with respect to the Euclidean variable on the horocycle space.
%%%%%%  Definition of Isometries  %%%%%%%%%%%%
\begin{Def}
 \label{Def:IsT}
For $s\in W$, define 
$T_{s}:L^{2}(X)\rightarrow L^{2}(\mathfrak{a}\times B,dHdb)$
by
$$
T_{s}u(H,b)=\int_{s\mathfrak{a}^{*}_{+}}
                  e^{i\lambda (H)}\mathcal{F}u(\lambda ,b)
                  \boldsymbol{c}^{-1}(\lambda)d\lambda
$$
and define
$T : L^{2}(X)\rightarrow L^{2}(\mathfrak{a}\times B,w^{-1}dHdb) $
by
$$
Tu(H,b)=\int_{\mathfrak{a}^{*}}
             e^{i\lambda (H)}\mathcal{F}u(\lambda ,b)
             \boldsymbol{c}^{-1}(\lambda)d\lambda .
$$
\end{Def}
%%%%%%%%%
%%%%%%%%%% Proposition for Isometry of T,Ts %%%%%%%%%%%
\begin{Prop}
 \label{Prop:IsT}
Both $T_{s}$ and $T$ are linear isometry operators respectively.
\end{Prop}
%%%%%%%%%%
%%%%%%%%%%% Proof of isometry %%%%%%%%%%%%%%%%%%%%%
\begin{proof}
It follows immediately from the $W$-invariantness of 
Fourier images of $L^{2}$ functions and the Plancherel theorem.
\end{proof}
%%%%%%%%%%%
%%%%%%%%%%%%%%%%%%%%%%%%%%%%%%%%%%%%%%%%%%%%%%%%%%%%
%%              proof of Theorem 1.2.             %%
%%%%%%%%%%%%%%%%%%%%%%%%%%%%%%%%%%%%%%%%%%%%%%%%%%%%
\noindent\textit{Proof of Theorem \ref{Thm:SmEf}}\,
%%%%%%%%%%%%%%   proof of (i)  %%%%%%%%%%%%%%%%%%%%%%%%%%%%%%%
(i) For any $s\in W$, we have
\begin{align*}
  T_{s}(p(D_{x})e^{ita(D_{x})}\psi)(H,b)
&=\int_{s\mathfrak{a}^{*}_{+}}
      e^{i\lambda (H)}
      \mathcal{F}(p(D_{x})e^{ita(D_{x})}\psi)(\lambda)
      \boldsymbol{c}^{-1}(\lambda)d\lambda \\
&=\int_{s\mathfrak{a}^{*}_{+}}
      e^{i\lambda (H)}
      p(s\lambda)e^{ita(s\lambda)}
      \mathcal{F}\psi(\lambda,b)
      \boldsymbol{c}^{-1}(\lambda)
  d\lambda \\
&=p(sD_{H})e^{ita(sD_{H})}
  (T_{s}\psi)(H,b).
\end{align*}
After decomposing $\mathfrak{a}^{*}$ as the
union of $s\mathfrak{a}^{*}_{+}(s=\pm 1)$,
we see that
\begin{align*}
 &\langle H\rangle^{-\delta}
  T( p(D_{x})e^{ita(D_{x})}\psi)(H,b)\\
=&\sum_{s\in W}
      \langle H\rangle^{-\delta}
      p(sD_{H})e^{ita(sD_{H})}(T_{s}\psi)(H,b).
\end{align*}
Then by taking $L^{2}$-norm over
$\mathbb{R}\times\mathfrak{a}\times B$
and applying the assumption, we obtain
\begin{align*}
     &\left\| p(D_{x})e^{ita(D_{x})}
             \psi
      \right\|_{L^{2}(\mathbb{R};L^{2,-\delta}(X) )}\\
 &\leq\sum_{s\in W}
      \left\| \langle H\rangle^{-\delta}
              p(sD_{H})e^{ita(sD_{H})}
              (T_{s}\psi)(H,b)
      \right\|_{L^{2}(\mathbb{R}\times\mathfrak{a}\times B,
                      w^{-1}dtdHdb)}\\
 &\leq w^{1/2}C_{\delta}\|\psi\|_{L^{2}(X)}.
\end{align*}
%%%%%%%%%%%%   proof of  (ii)  %%%%%%%%%%%%%%%%%%%%%%%%%%%%%%%
(ii) For $\chi$, we can take a cut off function
$\chi_{1}\in C^{\infty}(\mathfrak{a})$ such taht
$\chi_{1}(H)=0\,(|H|\leq 1/2)$,
$\chi_{1}=1$ on $\operatorname{supp}\chi$,
then $\chi=\chi\chi_{1}$ and we have
\begin{align*}
 &T
  \left(
   \chi(D_{x})q(D_{x})
   \int_{0}^{t}
     e^{i(t-\tau)a(D_{x})}
     f(\tau)
   d\tau
  \right)(H,b)\\
&=\sum_{s\in W}
  \chi(s D_{H})q(s D_{H})
   \int_{0}^{t}
     e^{i(t-\tau)a(D_{H})}
     T_{s}
     \left(
           \chi_{1}(D_{x})f_{\tau}
     \right)
     (H,b)
   d\tau.
\end{align*}
Then by applying the assumption, we obtain
\begin{align*}
 &\left\| \chi(D_{x})q(D_{x})
          \int_{0}^{t}
           e^{i(t-\tau)a(D_{H})}
           f_{\tau}
          d\tau
  \right\|_{L^{2}(\mathbb{R};L^{2,-\delta}(X))}\\
=&\sum_{s\in\{\pm 1\}}
  \left\| \langle H\rangle^{-\delta}
          \chi(s D_{H})q(s D_{H})
          \int_{0}^{t}
           e^{i(t-\tau)a(s D_{H})}
           T_{s}
           \left(
            \chi_{1}(D_{x})
            f_{\tau} \right)
           (H,b)
          d\tau
  \right\|_{L^{2}(\mathbb{R}\times\mathfrak{a}\times B,
                  w^{-1}dtdHdb)}\\
\leq &\sum_{s\in\{\pm 1\}}
      C_{\delta}
      \| \langle H\rangle^{\delta}
         T_{s}
         \left(
          \chi_{1}(D_{x})
          f_{t}
         \right)
      \|_{L^{2}(\mathbb{R}\times\mathfrak{a}\times B,
                w^{-1}dtdHdb)}.
\end{align*}
Here we put
$\chi_{2}(\lambda)
=\chi_{\mathfrak{a}_{+}^{*}}(\lambda)
 \chi_{1}(\lambda)
 \in C^{\infty}(\mathfrak{a}^{*})$, 
then we have 
$$
\langle H\rangle^{\delta}
 T_{s}
  \left( \chi_{1}(D_{x})f_{t}
   \right)
=\langle H\rangle^{\delta}
 \chi_{2}(D_{H})
 \langle H\rangle^{-\delta}
 \left( \langle H\rangle^{\delta}
        T(f_{t})
  \right).
$$
Since
$\langle H\rangle^{\delta}
 \chi_{2}(D_{H})
 \langle H\rangle^{-\delta}$
is an $L^{2}$-bounded operator on $\mathfrak{a}$, 
we obtain
\begin{align*}
&\left\| \chi(D_{x})q(D_{x})
          \int_{0}^{t}
           e^{i(t-\tau)a(D_{H})}
           f_{\tau}
          d\tau
  \right\|_{L^{2}(\mathbb{R};L^{2,-\delta}(X))}\\
&\leq C_{\delta ,\chi_{2}}
  \left\| \langle H\rangle^{\delta}
            T(f_{t})
  \right\|_{L^{2}(\mathbb{R}\times\mathfrak{a}\times B,
                  w^{-1}dtdHdb)}\\
&=C_{\delta,\chi_{2}}
  \left\| f
   \right\|_{L^{2}(\mathbb{R};L^{2,\delta}(X))}.
\end{align*}
So we complete the proof. 
\hfill $\square$
%%% end proof %%

%%%%%%%%%%%%%%%%%%%%%%%%%%%%%%%%%%%%%%%%%%%%%%%%%%%%%%%%%%%%%%%%%%%%%%
%  Section 5 : Gain of regularity for the Schroedinger   equation    %
%%%%%%%%%%%%%%%%%%%%%%%%%%%%%%%%%%%%%%%%%%%%%%%%%%%%%%%%%%%%%%%%%%%%%%
\section{Gain of regularity for the Schr\"{o}dinger evolution equation
         \label{GOR}}
In this section, we prove a gain of regularity for the Schr\"{o}dinger
evolution equation on real rank one symmetric spaces of noncompact type.
We also use the linear isometry introduced in the previous section
to reduce the argument to the Euclidean case.
We can obtain a gain of regularity with respect to the Euclidean variable
then we can recover the regularity of the solution
from its image of Radon transform by treating dual Radon transform
as an elliptic Fourier integral operator.

The following estimates on one-dimensional 
Euclidean spaces are well known.

%%%%%%%%%% Proposition for gain on one-dimensional %%%%%%%%%%%%%%%%%%%
\begin{Prop}
\label{Prop:GainConti}
On one-dimensional Euclidean space
$\mathbb{R}^{1}$,
for any $k\in\mathbb{N}$ 
and $\delta > 1/2$, 
we have the following continuous maps: 
\begin{alignat}{2}
 \langle x\rangle^{-k}
 L^{2}(\mathbb{R}^{1})
 \ni \phi 
&\mapsto
 t^{k}
 \langle {x}\rangle^{-k-\delta}
 \langle D_{x}\rangle^{k+1/2}
 e^{-it\Delta_{\mathbb{R}^{1}}}\phi
&\hspace{4pt}
 \in
&\hspace{4pt}
 L^{2}_{\mathrm{loc}}
 (\mathbb{R};
  L^{2}(\mathbb{R}^{1}),\label{eq:GainOnR1}\\
 \langle x\rangle^{-k}
 L^{2}(\mathbb{R}^{1})
 \ni \phi 
&\mapsto
 t^{k}
 \langle {x}\rangle^{-k\vphantom{-\delta}}
 \langle D_{x}\rangle^{k}
 e^{-it\Delta_{\mathbb{R}^{1}}}\phi
&\hspace{4pt} 
 \in
&\hspace{4pt}
 C(\mathbb{R};
 L^{2}(\mathbb{R}^{1})).
\end{alignat}
\end{Prop}
%%%%%%%%%%
We can also translate above estimates
on the one-dimensinal Euclidean space
into its on symmetric spaces by using 
isometry $T$.

\vspace{5pt}
%%%%%%%%%%%%%%%%%%%%%%%%%%%%%%%%%%%%%%%%%%%%%%%%%%%%%
%%           proof of Theorem 1.5                  %%
%%%%%%%%%%%%%%%%%%%%%%%%%%%%%%%%%%%%%%%%%%%%%%%%%%%%%
\noindent\textit{Proof of Theorem \ref{Thm:GRe}}
\, For any $\varphi\in L^{2}(X)$ with
$$\left\| \langle D_{\lambda}\rangle^{k}
    \big( \mathcal{F}\varphi(\lambda,b)
       \boldsymbol{c}^{-1}(\lambda)
    \big)
 \right\|_{L^{2}(\mathfrak{a}^{*}\times\Theta,d\lambda db)}
< \infty, $$ 
where $\Theta$ is an open subset of $B$, 
we can find 
\begin{align*}
T\big( e^{-it\Delta_{X}}\varphi\big)(H,b)
&=e^{\rho (H)}\Lambda\mathcal{R}
   \big( e^{-it\Delta_{X}}\varphi\big)(H,b)\\
&=e^{it|\rho|^{2}}e^{-it\Delta_{\mathfrak{a}}}\big( T\varphi\big)(H,b)
\end{align*}
and
\begin{align*}
 \left\| \langle H\rangle^{k}\big(T\varphi\big)(H,b)
  \right\|_{L^{2}(\mathfrak{a}\times\Theta,dHdb)}
&=\left\| \langle D_{\lambda}\rangle^{k}
   \big( \mathcal{F}\varphi(\lambda ,b)
         \boldsymbol{c}^{-1}(\lambda)
   \big)
  \right\|_{L^{2}(\mathfrak{a}^{*}\times\Theta,d\lambda db)}\\
&< \infty.
\end{align*}
Therefore by using 
Proposition \ref{Prop:GainConti} 
we have,
$$
\phi\in \langle x\rangle^{-k}
        L^{2}(\mathbb{R}^{1})
\Rightarrow
e^{-it\Delta_{\mathbb{R}^{1}}}
\phi
\in
H_{\mathrm{loc}}^{k+1/2}(\mathbb{R}^{1})
\quad \text{for\, a.e.\ }\; t\neq 0 ,
$$
hence we get 
\begin{equation}
\Lambda\mathcal{R}
\big( 
     e^{-it\Delta_{X}}\varphi
\big)
\in L^{2}( \Theta;H_{\mathrm{loc}}^{k+1/2}(\mathfrak{a}))
\quad \text{for\, a.e.\ }\; t\neq 0.
\label{eq:ReR}
\end{equation}
Now take $\chi_{1}\in\mathscr{D}(X)$ with 
$
\operatorname{supp}\chi_{1}
\subset 
\overline{B(o,2)}
$, 
$\chi_{1}=1$ 
near 
$B(o,1)$, and 
$\chi_{2}\in\mathscr{D}(\mathfrak{a})$ with 
$\chi_{2}=1$ 
near 
$\{H\in\mathfrak{a};|H|\leq 1\}$, 
then we can write
$$
\chi_{1}e^{-it\Delta_{X}}\varphi
=\chi_{1}\mathcal{R}^{*}\chi_{2}\Bar\Lambda
 \big(\Lambda\mathcal{R}(e^{-it\Delta_{X}}\varphi)\big).
$$
For $(o;\omega_{b}(o))\in T_{o}^{*} X\setminus 0$
with $b\in\Theta\cap (-\mathrm{id}_{B})\Theta$,
by using Theorem \ref{Thm:SoWF} we have
\begin{align*}
&(o;\pm\omega_{b}(o))
 \notin 
 \Gamma^{t}
 \circ
 \operatorname{WF}^{k+1/2}
 \big( 
       \Lambda\mathcal{R}(e^{-it\Delta_{X}}\varphi)
 \big)\\
\Rightarrow
&(o;\pm\omega_{b}(o))
 \notin 
 \Gamma^{t}
 \circ
 \operatorname{WF}^{k+1/2-\dim N/2}
 \big( \Bar\Lambda
 \big(\Lambda\mathcal{R}(e^{-it\Delta_{X}}\varphi)\big)\big)\\
\Rightarrow
&(o;\pm\omega_{b}(o))
 \notin 
 \operatorname{WF}^{k+1/2}
 \big( \mathcal{R}^{*}\chi_{2}\Bar\Lambda
 \big(\Lambda\mathcal{R}(e^{-it\Delta_{X}}\varphi)\big)\big)\\
\Leftrightarrow
&(o;\pm\omega_{b}(o))
 \notin 
 \operatorname{WF}^{k+1/2}
 \big(e^{-it\Delta_{X}}\varphi\big).
\end{align*}
Since
$$
 \Gamma^{t}
 \cap
 T^{*}_{(o,b,0)}(\Xi\times X)
=\left\{ (o,b,0;\lambda ( \omega_{b}(o)+dH))
    ; \lambda\in\mathbb{R}\setminus\{0\}
 \right\}, 
$$
we obtain 
\begin{align*}
&(o;\pm\omega_{b}(o))
 \notin 
 \Gamma^{t}
 \circ 
 \operatorname{WF}^{k+1/2}
 \big(\Lambda\mathcal{R}(e^{-it\Delta_{X}}\varphi)\big)\\
\Leftrightarrow
&(o,\pm b;dH)
 \notin 
 \operatorname{WF}^{k+1/2}
 \big(\Lambda\mathcal{R}(e^{-it\Delta_{X}}\varphi)\big).
\end{align*}
So from (\ref{eq:ReR}), we get
$
(o;\pm \omega_{b}(o))
 \notin
 \operatorname{WF}^{k+1/2}
\big(e^{-it\Delta_{X}}\varphi\big)
$
for a.e.\ $t\neq 0$.

Next, by using $G$-invariantness 
of the Schr\"{o}dinger evolution group $e^{-it\Delta_{X}}$
we show at other points. 
For any $g\in G$, we have
$$
\mathcal{F}(\tau_{g}\varphi)(\lambda,b)
=e^{(i\lambda -\rho)(A(g\cdot o,g\cdot b))}
 \mathcal{F}\varphi(\lambda,g\cdot b),
$$
where $\tau_{g}\varphi(x)=\varphi(g\cdot x)$.
For the non-Euclidean metric $A(x,b)$
we have the following formula
for the $G$-action:
$$
 A(g\cdot x,g\cdot b)
=A(x,b)+A(g\cdot o,g\cdot b).
$$
Also for the action on $B=G/MAN$ of $G$, 
we have 
$d(g\cdot b)
 =e^{-2\rho (A(g\cdot o,g\cdot b))}db$.
Then for $\tau_{g}\varphi$, we see that 
\begin{align*}
&\left\|
     \langle D_{\lambda}\rangle^{k}
     \big( 
        \mathcal{F}(\tau_{g}\varphi)(\lambda,b)
        \boldsymbol{c}^{-1}(\lambda)
     \big)
 \right\|_{L^{2}(\mathfrak{a}^{*}
                 \times g^{-1}
                 \cdot\Theta,
                 d\lambda db)}\\
=&\left\|
      \langle D_{\lambda}\rangle^{k}
      \big(
         e^{(i\lambda -\rho)(A(g\cdot o,g\cdot b))}
         \mathcal{F}\varphi(\lambda,g\cdot b)
         \boldsymbol{c}^{-1}(\lambda) 
      \big)
  \right\|_{L^{2}(\mathfrak{a}^{*}
                  \times g^{-1}
                  \cdot\Theta,
                  d\lambda db)}\\
=&\left\|
      \langle D_{\lambda}\rangle^{k}
      \big(
         e^{(i\lambda -\rho)(A(g\cdot o,b))}
         \mathcal{F}\varphi(\lambda,b)
         \boldsymbol{c}^{-1}(\lambda)
      \big) 
      e^{-\rho(A(g^{-1}\cdot o,g^{-1}\cdot b))}
  \right\|_{L^{2}(\mathfrak{a}^{*}
                  \times\Theta,
                  d\lambda db)},\\
\intertext{ by using
            $A(g^{-1}\cdot o,g^{-1}\cdot b)
            =-A(g\cdot o,b)$, we get
          }
=&\left\|
      \langle D_{\lambda}\rangle^{k}
      \big(e^{i\lambda (A(g\cdot o,b))}
         \mathcal{F}\varphi(\lambda,b)
         \boldsymbol{c}^{-1}(\lambda) 
      \big) 
  \right\|_{L^{2}(\mathfrak{a}^{*}\times\Theta ,d\lambda db)}\\
=&\left\|
      \big(
         \langle D_{\lambda}\rangle^{k}
         e^{i\lambda (A(g\cdot o,b))} 
         \langle D_{\lambda}\rangle^{-k} 
      \big)
      \langle D_{\lambda}\rangle^{k}
      \big(
         \mathcal{F}\varphi(\lambda ,b)
         \boldsymbol{c}^{-1}(\lambda)
      \big)
  \right\|_{L^{2}(\mathfrak{a}^{*}\times\Theta ,d\lambda db)}.
\end{align*}
Since 
$\{ \langle D_{\lambda}\rangle^{k}
    e^{i\lambda (A(g\cdot o,b))}
    \langle D_{\lambda}\rangle^{-k}
 \}_{b\in B} \in\mathcal{L}(L^{2}(\mathfrak{a}^{*}))$ 
is a bounded subset, so we get 
$$
\left\|
   \langle D_{\lambda}\rangle^{k}
   \big(
       \mathcal{F}(\tau_{g}\varphi)(\lambda ,b)
       \boldsymbol{c}^{-1}(\lambda )
   \big)
\right\|_{L^{2}(\mathfrak{a}^{*}\times 
                 g^{-1}\Theta ,d\lambda db)}
< \infty .
$$
Hence $\tau_{g}\varphi\in L^{2}(X)$
and $g^{-1}\Theta\subset B$ 
satisfies the assumption (\ref{eq:AsGR}).
By applying the above argument at the origin $o$, 
for any $b\in g^{-1}\Theta\cap (-\mathrm{id}_{B})g^{-1}\Theta$, 
and a.e.\ $t\neq 0$ we obtain 
\begin{align*}
(o;\pm\omega_{b}(o))
&\notin \operatorname{WF}^{k+1/2}
        \big(e^{-it\Delta_{X}}(\tau_{g}\varphi)\big)\\
&=      \operatorname{WF}^{k+1/2}
        \big(\tau_{g}(e^{-it\Delta_{X}}\varphi)\big)\\
&=      (\tau_{g})^{*}
        \big(
              \operatorname{WF}^{k+1/2}
              (e^{-it\Delta_{X}}\varphi)
        \big),
\end{align*}
hence 
\begin{align*}
(g\cdot o;\pm\omega_{g\cdot b}(g\cdot o))
&=(\tau_{g^{-1}})^{*}
 \big( (o;\omega_{b}(o))\big)\\
&\notin \operatorname{WF}^{k+1/2}(e^{-it\Delta_{X}}\varphi).
\end{align*}
Finally, if $\varphi\in L^{2,\delta}(X)$, 
then we can take $\Theta=B$, so the assertion folds. 
For gain of $H_{\mathrm{loc}}^{k}$-regularity, 
we can also prove by the same argument.
So we complete the proof.
\hfill $\square$
%%% end of proof %%%
\newpage

Finally, we prove Theorem \ref{Thm:GainConti}.
\vspace{5pt}

%%%%%%%%%%%%%%%%%%%%%%%%%%%%%%%%%
%%% Proof of Theorem 1.6      %%%
%%%%%%%%%%%%%%%%%%%%%%%%%%%%%%%%%
\noindent\textit{Proof of Theorem \ref{Thm:GainConti}}
\, For any $\varphi\in L^{2,k}(X)$, from the definition 
of the weighted $L^{2}$-norm, we have
\begin{align*}
 &\|
    t^{k}
    \langle D_{x}\rangle^{k+1/2}
    e^{-it\Delta_{X}}
    \varphi
  \|_{L^{2,-k-\delta}(X)}\\
=&\| 
    t^{k}
    \langle H\rangle^{-k-\delta}
    T
    \left( 
         \langle D_{x}\rangle^{k+1/2}
         e^{-it\Delta_{X}}
         \varphi
    \right)
  \|_{L^{2}(\mathfrak{a}
            \times B,
            w^{-1}dHdb)}\\
=&\| 
    t^{k}
    \langle H\rangle^{-k-\delta}
    \langle D_{H}\rangle^{k+1/2}
    e^{-it\Delta_{\mathfrak{a}}+it|\rho|^{2}}
    T\varphi
  \|_{L^{2}(\mathfrak{a}
            \times B,
            w^{-1}dHdb)}.
\end{align*}
For any $T>0$, by using the continuity
of the map (\ref{eq:GainOnR1}), we obtain 
\begin{align*}
 &\| 
    t^{k}
    \langle D_{x}\rangle^{k+1/2}
    e^{-it\Delta_{X}}
    \varphi
  \|_{L^{2}((-T,T)
      ;L^{2,-k-\delta}(X))}\\
=&\| 
    t^{k} 
    \langle H\rangle^{-k-\delta}
    \langle D_{H}\rangle^{k+1/2}
    e^{-it\Delta_{\mathfrak{a}}}
    (T\varphi)
  \|_{L^{2}((-T,T)       
            \times
            \mathfrak{a} 
            \times
            B;w^{-1}dHdb
           )}\\
\leq
  &C_{T}
  \| 
     \langle H\rangle^{k}
     (T\varphi)
  \|_{L^{2}(\mathfrak{a}
            \times B;w^{-1}dHdb)}\\
=&C_{T}
  \|\varphi\|_{L^{2,k}(X)},
\end{align*}
for some $C_{T}>0$ 
depending on $T$.
The rest assertion
proved by the same argument.
\hfill $\square$
%%% end of proof %%%
%%%%%%%%%%%%%%%%%%%%%%%%%%%%%%%%%
%%%         References        %%%
%%%%%%%%%%%%%%%%%%%%%%%%%%%%%%%%%

%%%%%%%%%%%%%%%%%%%%%%%%%
% address & E-Mail      %
%%%%%%%%%%%%%%%%%%%%%%%%%
\address{
Mathematical Institute\\
Tohoku University\\
Sendai 980-8578\\
Japan
}
 {sa6m13@math.tohoku.ac.jp}
%%%%%%%%%%%%%%
%%%%%%%%%%%%%%
%%%%%%%%%%%%%%

\begin{bibdiv}
\begin{biblist}
\bib{BenKla}{article}{
   author={Ben-Artzi, Matania},
   author={Klainerman, Sergiu},
   title={Decay and regularity for the Schr\"odinger equation},
  %note={Festschrift on the occasion of the 70th birthday of Shmuel Agmon},
   journal={J. Anal. Math.},
   volume={58},
   date={1992},
   pages={25--37},
  %issn={0021-7670},
  %review={\MR{1226935 (94e:35053)}},
}
\bib{Beylkin}{article}{
   author={Beylkin, Gregory},
   title={The inversion problem and applications of the generalized Radon
   transform},
   journal={Comm. Pure Appl. Math.},
   volume={37},
   date={1984},
   number={5},
   pages={579--599},
  %issn={0010-3640},
  %review={\MR{752592 (86a:44002)}},
}
\bib{Chihara}{article}{
   author={Chihara, Hiroyuki},
   title={Smoothing effects of dispersive pseudodifferential equations},
   journal={Comm. Partial Differential Equations},
   volume={27},
   date={2002},
   number={9-10},
   pages={1953--2005},
  %issn={0360-5302},
  %review={\MR{1941663 (2004c:35441)}},
}
\bib{CKS}{article}{
   author={Craig, Walter},
   author={Kappeler, Thomas},
   author={Strauss, Walter},
   title={Microlocal dispersive smoothing for the Schr\"odinger equation},
   journal={Comm. Pure Appl. Math.},
   volume={48},
   date={1995},
   number={8},
   pages={769--860},
  %issn={0010-3640},
  %review={\MR{1361016 (96m:35057)}},
}
\bib{Doi0}{article}{
   author={Doi, Shin-ichi},
   title={Remarks on the Cauchy problem for Schr\"odinger-type equations},
   journal={Comm. Partial Differential Equations},
   volume={21},
   date={1996},
   number={1-2},
   pages={163--178},
  %issn={0360-5302},
  %review={\MR{1373768 (96m:35058)}},
}
\bib{Doi1}{article}{
   author={Doi, Shin-ichi},
   title={Smoothing effects of Schr\"odinger evolution groups on Riemannian
   manifolds},
   journal={Duke Math. J.},
   volume={82},
   date={1996},
   number={3},
   pages={679--706},
  %issn={0012-7094},
  %review={\MR{1387689 (97f:58141)}},
}
\bib{Doi2}{article}{
   author={Doi, Shin-ichi},
   title={Commutator algebra and abstract smoothing effect},
   journal={J. Funct. Anal.},
   volume={168},
   date={1999},
   number={2},
   pages={428--469},
   issn={0022-1236},
  %review={\MR{1719229 (2001e:58036)}},
}
\bib{Doi3}{article}{
   author={Doi, Shin-ichi},
   title={Smoothing effects for Schr\"odinger evolution equation and global
   behavior of geodesic flow},
   journal={Math. Ann.},
   volume={318},
   date={2000},
   number={2},
   pages={355--389},
  %issn={0025-5831},
  %review={\MR{1795567 (2001h:58045)}},
}
\bib{EgOk}{article}{
   author={Eguchi, M.},
   author={Okamoto, K.},
   title={The Fourier transform of the Schwartz space on a symmetric space},
   journal={Proc. Japan Acad. Ser. A Math. Sci.},
   volume={53},
   date={1977},
   number={7},
   pages={237--241},
  %issn={0386-2194},
  %review={\MR{0499949 (58 \#17691)}},
}
\bib{GaVa}{book}{
   author={Gangolli, Ramesh},
   author={Varadarajan, V. S.},
   title={\rm{``Harmonic analysis of spherical functions 
                on real reductive groups''}},
   series={Ergebnisse der Mathematik und ihrer Grenzgebiete [Results in
   Mathematics and Related Areas]},
   volume={101},
   publisher={Springer-Verlag},
   place={Berlin},
   date={1988},
   pages={xiv+365},
  %isbn={3-540-18302-7},
  %review={\MR{954385 (89m:22015)}},
}
\bib{GonzalezQuinto}{article}{
   author={Gonzalez, Fulton},
   author={Quinto, Eric Todd},
   title={Support theorems for Radon transforms on higher rank symmetric
   spaces},
   journal={Proc. Amer. Math. Soc.},
   volume={122},
   date={1994},
   number={4},
   pages={1045--1052},
  %issn={0002-9939},
  %review={\MR{1205492 (95b:44002)}},
}
\bib{GuilSter}{book}{
   author={Guillemin, Victor},
   author={Sternberg, Shlomo},
   title={\rm{``Geometric asymptotics''}},
   note={Mathematical Surveys, No. 14},
   publisher={American Mathematical Society},
   place={Providence, R.I.},
   date={1977},
   pages={xviii+474 pp. (one plate)},
  %review={\MR{0516965 (58 \#24404)}},
}
\bib{GreUhl}{article}{
   author={Greenleaf, Allan},
   author={Uhlmann, Gunther},
   title={Microlocal techniques in integral geometry},
   conference={
      title={Integral geometry and tomography},
      address={Arcata, CA},
      date={1989},
   },
   book={
      series={Contemp. Math.},
      volume={113},
      publisher={Amer. Math. Soc.},
      place={Providence, RI},
   },
   date={1990},
   pages={121--135},
  %review={\MR{1108649 (93d:44003)}},
}
\bib{HaCh}{article}{
   author={Harish-Chandra},
   title={Discrete series for semisimple Lie groups. II. Explicit
   determination of the characters},
   journal={Acta Math.},
   volume={116},
   date={1966},
   pages={1--111},
  %issn={0001-5962},
  %review={\MR{0219666 (36 \#2745)}},
}
\bib{Helgason1}{article}{
   author={Helgason, Sigurdur},
   title={Fundamental solutions of invariant differential operators on
   symmetric spaces},
   journal={Amer. J. Math.},
   volume={86},
   date={1964},
   pages={565--601},
  %issn={0002-9327},
  %review={\MR{0165032 (29 \#2323)}},
}

\bib{Helgason2}{article}{
   author={Helgason, Sigurdur},
   title={A duality for symmetric spaces with applications to group
   representations},
   journal={Advances in Math.},
   volume={5},
   date={1970},
   pages={1--154},
  %issn={0001-8708},
  %review={\MR{0263988 (41 \#8587)}},
}
\bib{Helgason3}{article}{
   author={Helgason, Sigurdur},
   title={Paley-Wiener theorems and surjectivity of invariant differential
   operators on symmetric spaces and Lie groups},
   journal={Bull. Amer. Math. Soc.},
   volume={79},
   date={1973},
   pages={129--132},
  %review={\MR{0312158 (47 \#720)}},
}
\bib{HeGGA}{book}{
   author={Helgason, Sigurdur},
   title={\rm{``Groups and geometric analysis''}},
   series={Pure and Applied Mathematics},
   volume={113},
   note={Integral geometry, invariant differential operators, and spherical
   functions},
   publisher={Academic Press Inc.},
   place={Orlando, FL},
   date={1984},
   pages={xix+654},
  %isbn={0-12-338301-3},
  %review={\MR{754767 (86c:22017)}},
}
\bib{HeGASS}{book}{
   author={Helgason, Sigurdur},
   title={\rm{``Geometric analysis on symmetric spaces''}},
   series={Mathematical Surveys and Monographs},
   volume={39},
   publisher={American Mathematical Society},
   place={Providence, RI},
   date={1994},
   pages={xiv+611},
  %isbn={0-8218-1538-5},
  %review={},
}
\bib{HeDLS}{book}{
   author={Helgason, Sigurdur},
   title={\rm{``Differential geometry, Lie groups, and symmetric spaces''}},
   series={Graduate Studies in Mathematics},
   volume={34},
   note={Corrected reprint of the 1978 original},
   publisher={American Mathematical Society},
   place={Providence, RI},
   date={2001},
   pages={xxvi+641},
  %isbn={0-8218-2848-7},
  %review={\MR{1834454 (2002b:53081)}},
}
\bib{Ho-III}{book}{
   author={H{\"o}rmander, Lars},
   title={\rm{``The analysis of linear partial differential operators. III''}},
  %series={Grundlehren der Mathematischen Wissenschaften [Fundamental
  %Principles of Mathematical Sciences]},
   volume={274},
  %note={Pseudodifferential operators},%
   publisher={Springer-Verlag},
   place={Berlin},
   date={1985},
   pages={viii+525},
  %isbn={3-540-13828-5},
  %review={},
}
\bib{Hoshiro}{article}{
   author={Hoshiro, Toshihiko},
   title={Decay and regularity for dispersive equations with constant
   coefficients},
   journal={J. Anal. Math.},
   volume={91},
   date={2003},
   pages={211--230},
  %issn={0021-7670},
  %review={\MR{2037408 (2005a:35005)}},
}
\bib{Kato}{article}{
   author={Kato, Tosio},
   title={On the Cauchy problem for the (generalized) Korteweg-de Vries
   equation},
   conference={
      title={Studies in applied mathematics},
   },
   book={
      series={Adv. Math. Suppl. Stud.},
      volume={8},
      publisher={Academic Press},
      place={New York},
   },
   date={1983},
   pages={93--128},
  %review={\MR{759907 (86f:35160)}},
}
\bib{KatoYajima}{article}{
   author={Kato, Tosio},
   author={Yajima, Kenji},
   title={Some examples of smooth operators and the associated smoothing
   effect},
   journal={Rev. Math. Phys.},
   volume={1},
   date={1989},
   number={4},
   pages={481--496},
  %issn={0129-055X},
  %review={\MR{1061120 (91i:47013)}},
}
\bib{Kumano-go}{book}{
   author={Kumano-go, H.},
   title={\rm{``Pseudodifferential operators''}},
   note={},
   publisher={MIT Press},
   place={Cambridge, Mass.},
   date={1981},
   pages={xviii+455},
  %isbn={0-262-11080-6},
  %review={},
}
\bib{QuintoRadon}{article}{
   author={Quinto, Eric Todd},
   title={Real analytic Radon transforms on rank one symmetric spaces},
   journal={Proc. Amer. Math. Soc.},
   volume={117},
   date={1993},
   number={1},
   pages={179--186},
  %issn={0002-9939},
  %review={\MR{1135080 (93e:44005)}},
}
\bib{QuintoXray}{article}{
   author={Quinto, Eric Todd},
   title={An introduction to X-ray tomography and Radon transforms},
   conference={
      title={The Radon transform, inverse problems, and tomography},
   },
   book={
      series={Proc. Sympos. Appl. Math.},
      volume={63},
      publisher={Amer. Math. Soc.},
      place={Providence, RI},
   },
   date={2006},
   pages={1--23},
  %review={\MR{2208234 (2007i:35232)}},
}
\bib{RodTao}{article}{
   author={Rodnianski, I.},
   author={Tao, T.},
   title={Longtime decay estimates for the Schr\"odinger equation on
   manifolds},
   conference={
      title={Mathematical aspects of nonlinear dispersive equations},
   },
   book={
      series={Ann. of Math. Stud.},
      volume={163},
      publisher={Princeton Univ. Press},
      place={Princeton, NJ},
   },
   date={2007},
   pages={223--253},
  %review={\MR{2333213 (2008g:58035)}},
}
\bib{Shubin}{book}{
   author={Shubin, M. A.},
   title={\rm{``Pseudodifferential operators and spectral theory''}},
   edition={2},
  %note={Translated from the 1978 Russian original by Stig I. Andersson},
   publisher={Springer-Verlag},
   place={Berlin},
   date={2001},
   pages={xii+288},
  %isbn={3-540-41195-X},
  %review={\MR{1852334 (2002d:47073)}},
}
\bib{Sugimoto}{article}{
   author={Sugimoto, Mitsuru},
   title={Global smoothing properties of generalized Schr\"odinger
   equations},
   journal={J. Anal. Math.},
   volume={76},
   date={1998},
   pages={191--204},
  %issn={0021-7670},
  %review={\MR{1676995 (2000a:35033)}},
}
\bib{Tits}{article}{
   author={Tits, J.},
   title={Sur certaines classes d'espaces homog\`enes de groupes de Lie},
   language={French},
   journal={Acad. Roy. Belg. Cl. Sci. M\'em. Coll. in $8\sp \circ$},
   volume={29},
   date={1955},
   number={3},
   pages={268},
  %review={\MR{0076286 (17,874f)}},
}
 \end{biblist}
 \end{bibdiv}
\end{document}